\documentclass[11pt,a4paper]{article}
\usepackage[T1]{fontenc}
\usepackage[utf8]{inputenc}
\usepackage{lmodern,amsmath,amssymb,amsthm,mathtools}
\usepackage[a4paper,margin=29mm]{geometry}
\usepackage{microtype,booktabs,array}
\usepackage{scalerel}
\usepackage{pgfplots}
\pgfplotsset{compat=1.17}
\usepackage[noautomatic]{imakeidx}
\makeindex[intoc,columns=2]
\usepackage[hidelinks]{hyperref}
\allowdisplaybreaks[2]
\newtheorem{theorem}{Theorem}[section]
\newtheorem{lemma}[theorem]{Lemma}
\newtheorem{proposition}[theorem]{Proposition}
\newtheorem{corollary}[theorem]{Corollary}
\theoremstyle{definition}
\newtheorem{definition}[theorem]{Definition}
\theoremstyle{remark}
\newtheorem{remark}[theorem]{Remark}
\newcommand{\N}{\mathbb{N}}
\newcommand{\Z}{\mathbb{Z}}
\newcommand{\Q}{\mathbb{Q}}
\newcommand{\hh}{\widehat{h}}
\newcommand{\ord}{\mathcal{C}}
\newcommand{\Pow}{\operatorname{Pow}_2}
\newcommand{\Sol}{\operatorname{Sol}}
\newcommand{\va}{\mathbf{a}}
\newcommand{\vw}{\mathbf{w}}
\newcommand{\vz}{\mathbf{z}}

\def\sqdot{\mathbin{\scalerel*{\strut\rule{1.5ex}{1.5ex}}{\cdot}}}
\newcommand{\OO}{\mathcal{O}}
\newcommand{\GG}{\mathbf{G}}

\usepackage{etoolbox}
\makeatletter
\patchcmd{\l@section}{\addvspace{1.0em \@plus\p@}}{\addvspace{0.45em \@plus\p@}}{}{}
\patchcmd{\l@part}{\addvspace{2.25em \@plus\p@}}{\addvspace{1.2em \@plus\p@}}{}{}
\makeatother

\title{On the proof of the Diophantine singlefold-ness conjecture\\[3pt]\large Canonical Height Quotients and Exact Witness Transfer}
\author{Domenico Cantone\\[2pt]
{\small Department of Mathematics and Computer Science (DMI)}\\
{\small University of Catania, Italy}\\
{\small\href{mailto:domenico.cantone@unict.it}{\texttt{domenico.cantone@unict.it}}}
\and
Eugenio G. Omodeo\\[2pt]
{\small Department of Mathematics, Informatics and Geosciences (DMG)}\\
{\small University of Trieste, Italy}\\
{\small\href{mailto:eomodeo@units.it}{\texttt{eomodeo@units.it}}}}
\date{November 27, 2026}
\begin{document}
\maketitle
\begin{abstract}
\noindent We give a modular simplification of the proof of the singlefold representation theorem established in OpenAI's September 2026 preprint. The elliptic-height method, applied to the torsion-free Mordell curve $y^2=x^3-74$, yields explicit arithmetic parameters without height amplification. A classical Pell quotient and Matiyasevich's one-exponential normal form replace the encoding of computation histories. We then develop two relative constructions, independent of the general singlefold theorem. First, for finitely generated Abelian groups of positive rank with singlefold coded operations and a height whose logarithm differs from a positive definite quadratic form by a bounded amount, we construct canonical cyclic quotients by finite comparisons. In particular, the height-window method applies to elliptic curves of arbitrary positive rank. Second, under suitable growth bounds we obtain a power specification whose complete witness fibre is in bijection with the admissible tuples of a supplied growth system. One shared cutoff suffices for a finite conjunction of power equations, preserving the bound on the multiplicity of the supplied fibre even for nonfunctional families and allowing sparse domains under additional hypotheses. These constructions address the two objectives proposed in our earlier work with Cuzziol (2024).
\end{abstract}

\newpage
{\small\setcounter{tocdepth}{2}\tableofcontents}
\clearpage

\part{A simpler proof of the singlefold conjecture}\label{part:proof}
\section{Introduction}
Hilbert's tenth problem\index{Hilbert's tenth problem} asked for a decision procedure for the solvability of polynomial equations over the integers. The negative answer came as a consequence of a theorem of considerably wider scope. According to the Davis--Putnam--Robinson--Matiyasevich theorem\index{DPRM theorem} \cite{DPR,Mat1970}, recursive enumerability (listability) and polynomial Diophantine representability capture the same family of relations on the natural numbers. This theorem does not, however, tell us how many solutions a Diophantine equation specifying a relation may have once its parameters have been assigned values.

Throughout this paper, $\N=\{0,1,2,\ldots\}$ and $\N_{>0}=\N\setminus\{0\}$. Here and below, $\log$ denotes the natural logarithm, unless a different base is explicitly indicated. We adopt the convention $0^0=1$. Unless otherwise specified, the parameters and unknowns of polynomial specifications range over $\N$. A set $S\subseteq\N^n$ is termed \emph{recursively enumerable} (\emph{r.e.}), or listable\index{recursively enumerable (r.e.)}\index{listable|see{recursively enumerable (r.e.)}}, when some algorithm accepts precisely its members, possibly running indefinitely on nonmembers. In contrast to the case of a decidable set, such an algorithm need not give a negative answer after finitely many steps. By a Diophantine representation of $S$\index{Diophantine representation} we mean a valid biimplication of the form
\begin{equation}\label{eq:repr}
 \va\in S\quad\Longleftrightarrow\quad \exists\vw\in\N^m\ P(\va,\vw)=0,
\end{equation}
where $P$ is a polynomial with coefficients in $\Z$. We call $\va$ the \emph{parameters}\index{parameter} and $\vw$ the auxiliary variables, or \emph{witnesses}\index{witness}. For fixed $\va$, the set $\Sol_{\N}(P(\va,\sqdot))=\{\vw\in\N^m:P(\va,\vw)=0\}$ is the \emph{complete witness fibre}\index{complete witness fibre}; the placeholder $\sqdot$ marks the witness argument. A tuple is complete when it includes every auxiliary coordinate, including those introduced in implementing constituent operations. A representation is termed \emph{finitefold}\index{finitefold representation} if every such fibre is finite, and \emph{univocal}, or singlefold\index{singlefold (univocal) representation}\index{univocal|see{singlefold (univocal) representation}}, if every fibre contains at most one tuple. More generally, for a fixed positive integer $k$, a representation is termed \emph{$k$-bounded}\index{k-bounded representation@$k$-bounded representation} if every complete witness fibre has cardinality at most $k$. These properties concern the specified representation. If exponentiation with variable exponents is also admitted in the defining expressions, the representation is termed exponential Diophantine\index{exponential Diophantine representation}.

Consider, for instance, the specification of $a\ne b$ obtained by splitting the inequality into its two alternatives $a>b$ and $b>a$. If we introduce a separate slack variable\index{slack variable} (a natural-valued auxiliary variable converting an inequality into an equality) for each alternative, we get
\[
 (a-b-x-1)(b-a-y-1)=0.
\]
Whenever $a\ne b$, this equation has infinitely many solving pairs $(x,y)$: the slack of the unused alternative is unrestricted. We could avoid this inconvenience by using the same variable $x$ in both factors, namely $(a-b-x-1)(b-a-x-1)=0$; its sole possible value would then be $|a-b|-1$. If, instead, we retain the two variables, it suffices to fix the unused one at zero:
\[
 \big((a-b-x-1)^2+y\big)\big((b-a-y-1)^2+x\big)=0.
\]
The resulting specification is univocal. This second device will be useful for more elaborate disjunctions (cf. Lemma~\ref{lem:compiler}).

In 1974, Matiyasevich \cite{Mat1974} showed that the exponential Diophantine representation of an r.e. relation can be made univocal. Can one reconcile representability with univocity without resorting to exponentiation? One route to a positive answer, advocated by Davis \cite{Davis1968} in a related setting and pursued in \cite{CCFO2019,CCFO2021,CO2021,CCO2023,CCO2024}, involves special quaternary quartic equations. Eight such equations were proposed as possible routes to finitefold representability; the arithmetic of those equations and the conditional reductions leading to them are distinct from the argument to be developed here.

Behind those investigations lie two complementary issues (cf. \cite{CCFO2021}, \cite[Appendix A, p.~27]{CO2021}, and \cite[Sections 6 and 7]{CCO2024}). What requirements should be imposed upon a binary (or dyadic) relation $M\subseteq\N^2$ for it to supersede exponentiation in a Diophantine specification without jeopardizing finitefold-ness, or even univocity? And can one exhibit a relation meeting those requirements that has, in turn, a polynomial Diophantine representation of the desired kind? Here a distinction is essential: a numerical growth condition has no bearing by itself on the number of solving tuples, and finiteness for each fixed pair of parameters does not ensure finiteness when only the first parameter is fixed. We shall return to these issues in Section~\ref{sec:cco-programme}.

An unexpected way of tackling the second issue was opened up by the manuscript of September 2026 \cite{OpenAI2026}: the heights of rational points on a rank-one elliptic curve supply a univocal growth function, from which exponentiation can be recovered by elementary estimates on Pell solutions. OpenAI also reports a Lean formalization of the main theorem of that preprint \cite{OpenAI2026}, see \cite{OpenAILean242}, and community records report external checks of its formal statement \cite{Community242,Atlas242}; the formalization concerns the singlefold theorem itself, and the arguments of the present paper are not formally verified.

The decisive idea of using elliptic heights and a fixed multiplicative interval around $1$ to recognize the square of an index is due to \cite[Section 4, Proposition 4.2]{OpenAI2026}. Part~\ref{part:proof} adapts this idea to the torsion-free Mordell curve $y^2=x^3-74$, whose group structure and height bound we certify by elementary arguments without relying on elliptic-curve tables or numerical rank and height computations. This yields explicit parameters for the arithmetic module; the classical Pell quotient and Matiyasevich's one-exponential normal form simplify the conclusion of the proof. Part~\ref{part:extensions} develops two further directions: a general height principle, including canonical cyclic quotients of naturally coded groups with the stated height and singlefold arithmetic hypotheses (Section~\ref{sec:indices-principle}), and exact witness transfer\index{exact witness transfer} from supplied growth systems to power specifications, including a shared cutoff for finite conjunctions (Section~\ref{sec:transfer}). Together with the growth function of Part~\ref{part:proof}, these constructions, relative in the sense specified at the start of Part~\ref{part:extensions}, address the two objectives set out in \cite[Section 7]{CCO2024}; their relation to earlier work is discussed in Section~\ref{sec:prior}.

Let us recall the following theorem.
\begin{theorem}[OpenAI; simplified proof]\label{thm:main}\index{singlefold representation theorem}
For every recursively enumerable set $S\subseteq\N^n$, $n\geq1$, there exist an integer $m\geq1$ and a polynomial $P\in\Z[A_1,\ldots,A_n,W_1,\ldots,W_m]$ such that
\[
 \#\Sol_{\N}(P(\va,\sqdot))=\begin{cases}1,&\va\in S,\\0,&\va\notin S.\end{cases}
\]

\end{theorem}

The proof rests on the following three steps, whose roles it is useful to keep separate:
\nopagebreak[4]
\begin{center}
\begin{tabular}{@{}c@{}}
A rank-one elliptic curve and its heights\\[2pt]
$\Downarrow$\\[2pt]
A singlefold graph of a function $g$ with $\log g(u)\sim C(\log u)^2$, $C>0$\\[2pt]
$\Downarrow$\quad (Pell equations)\\[2pt]
A singlefold graph of $c=2^n$\\[2pt]
$\Downarrow$\quad (Matiyasevich's 1974 theorem)\\[2pt]
Singlefold polynomial representations of all r.e. sets.
\end{tabular}
\end{center}
Sections~\ref{sec:systems}--\ref{sec:conclusion} establish Theorem~\ref{thm:main} in the order outlined above. The requisite details about rational-point codes and an elementary derivation of the classical quotient are relegated to Appendices~\ref{app:codes} and \ref{app:quotient}, respectively; Appendix~\ref{app:descent} establishes the rank of the chosen curve. Logarithms, signed point indices and real-valued bounds will aid us in the proofs; none of these is an extra operation or unknown in the resulting polynomial.

The reader interested solely in the singlefold conjecture may stop at the end of Part~\ref{part:proof}. Part~\ref{part:extensions} develops the general height principle (Section~\ref{sec:indices-principle}), exact witness transfer (Section~\ref{sec:transfer}), and a consequence for quartic equations (Section~\ref{sec:consequences}). It concludes by revisiting the programme proposed in our earlier work (Section~\ref{sec:cco-programme}).

\section{Singlefold systems and a classical reduction}\label{sec:systems}
We will often specify a relation by a finite system of equations rather than by a sole equation. A relation $R\subseteq\N^d$ has a singlefold polynomial system if the conjunction of equations $F_i(\va,\vz)=0$ has exactly one complete auxiliary tuple $\vz$ when $\va\in R$, and none otherwise. The complete witness fibre is then the set of tuples satisfying every equation. Since on natural (or integer) tuples the conjunction amounts to $\sum_iF_i^2=0$, this choice of notation entails no loss of univocity and calls for no extra variables.

By a \emph{natural coding}\index{natural coding} of a set we mean an injection into $\N^s$ for a fixed positive integer $s$; a code is valid if it belongs to the image. By singlefold validity\index{singlefold validity} we mean that the predicate of being a valid code has a singlefold polynomial representation.

By a function with a singlefold polynomial graph (or, for brevity, a univocal function) we mean a function with a singlefold representation of its graph: fixing an input in its domain determines not only the output but also every auxiliary value. As the examples in \cite[Section 2]{CCO2024} illustrate, conjunction and composition can safely be used, whereas disjunction requires some care lest the auxiliary variables of an unused alternative remain undetermined. We collect the elementary rules (using the selector mechanism of \cite[Section 2, Lemma 2.2]{OpenAI2026}) that will repeatedly come into play.
\begin{lemma}[Singlefold bookkeeping]\label{lem:compiler}\index{singlefold bookkeeping}
\hspace{0pt}
\begin{enumerate}
\item A finite conjunction of singlefold systems, with fresh internal variables for each call, has unique internal witnesses whenever all its displayed arguments are fixed and the conjunction is true.
\item A finite composition of functions with singlefold graphs has a singlefold graph.
\item A finite union of singlefold relations that are disjoint on their common displayed arguments has a singlefold system. For functions defined by cases, guards that partition the input domain and determine a unique output in each case suffice.
\item Comparisons of integer polynomial expressions, maxima of natural numbers, lexicographic comparisons of natural tuples of fixed length, and division with remainder by a positive natural number can be expressed with unique witnesses.
\end{enumerate}
\end{lemma}
\begin{proof}
For the first two assertions, determine the auxiliary values of each constituent specification once its arguments are fixed; in a composition the outputs are determined in succession. For a disjoint union with $r$ cases, use natural selectors\index{selector} $e_i$ subject only to $\sum_{i=1}^r e_i=1$: precisely one of them equals 1 and all the others equal 0. Multiply every equation $F=0$ belonging to case $i$ by $e_i$, and require $(1-e_i)z=0$ for every variable $z$ local to that case. The active case is unique; all inactive local coordinates are zero. Shared arguments are not zeroed. This convention applies also to every internal variable of a nested call.

As for comparisons, $F\geq G$ can be specified by $F-G=s$, and $F>G$ by $F-G=1+s$, with one natural slack $s$. The cases $a\geq b$ and $a<b$ define $\max\{a,b\}$ without overlap, assigning ties to the first case. A lexicographic comparison is a disjoint union of the cases distinguished by the first unequal coordinate, with the appropriate strict inequality and equality of earlier coordinates; for a nonstrict comparison include the equality case. Euclidean division is described by $a=dq+r$, $r<d$, $d>0$, which fixes $q,r$ and their comparison slacks uniquely.
\end{proof}

The reason for concentrating on powers of two is a classical theorem of Matiyasevich, in the one-exponential form recalled in \cite[Theorem 5]{CCO2024} (cf. \cite{Mat1974}).
\begin{theorem}[Matiyasevich, 1974]\label{thm:mat}\index{Matiyasevich's one-exponential normal form}
For every recursively enumerable $S\subseteq\N^n$ there are a nonnegative integer $s$ and an integer polynomial $D$ such that
\begin{equation}\label{eq:mat}
 \va\in S\quad\Longleftrightarrow\quad
 \exists u,v\in\N\;\exists\mathbf{x}\in\N^s\quad
 \bigl(2^u=v\ \wedge\ D(\va,\mathbf{x},u,v)=0\bigr).
\end{equation}
For each fixed $\va\in S$, the complete tuple $(\mathbf{x},u,v)$ is unique.
\end{theorem}
To obtain this form of the theorem, start from $P(\mathbf{a},\mathbf{x})=4^y+y$ recorded by Matiyasevich in \cite[p.~80, Eq.~(2), of the original publication]{Mat2010}, with a unique complete tuple $(\mathbf{x},y)$ on true inputs. That source uses positive integers. Apply its theorem to the translated r.e. set $S^+=\{\mathbf{a}+\mathbf 1:\mathbf{a}\in S\}$ and replace each positive witness $x_i$ by $X_i+1$ and $y$ by $Y+1$, with $X_i,Y\in\N$; here $\mathbf 1$ denotes the vector of ones of the appropriate dimension. Introduce $u=2(Y+1)$ and $v=2^u$. The two polynomial equations
\[
 u=2(Y+1),\qquad P(\va+\mathbf 1,\mathbf{X}+\mathbf 1)-v-Y-1=0
\]
combine into one equation by summing their squares. This gives \eqref{eq:mat}, with $\mathbf{X},Y$ among its witnesses; the translations and the additional output coordinates preserve uniqueness.

Any uniform bound on the number of complete witnesses in a polynomial representation of $2^u=v$ transfers to all recursively enumerable relations by substitution in \eqref{eq:mat}: the classical tuple is unique on true inputs, and each witness of the substituted system extends it exactly once. In particular, finitefold representability transfers in the same way. To prove Theorem~\ref{thm:main}, it therefore suffices to give a univocal polynomial specification of $2^u=v$. For the development leading to the classical reduction---a univocal version of Davis's normal form, the elimination of bounded universal quantification, and the reduction to one power of two---see \cite[Lemmas 4.1 and 5.3, Theorem 3 and Section 6]{CCO2024}.

\section{What elliptic heights can do}\label{sec:height}
\subsection{Recognizing the square of an index}
We first recall the index-selection argument of \cite[Section 4.2, Proposition 4.2]{OpenAI2026} in terms of an abstract sequence. Suppose we have a sequence $W_0,W_1,\ldots$ of positive integers satisfying
\begin{equation}\label{eq:model}
 \log W_j=aj^2+\varepsilon_j,\qquad |\varepsilon_j|\leq c,
\end{equation}
where $a>0$ and $c$ is independent of $j$. Suppose also that, once a term has been selected, we can select its neighbours and the term at twice its index without having to specify the index itself.

To exploit the bounded error, consider the identities
\[
 (2k)^2-2k^2=2k^2,\qquad (r+1)^2-r^2-1=2r.
\]
Subtracting the second from the first, and replacing the squares by the corresponding logarithms in \eqref{eq:model}, we obtain
\begin{equation}\label{eq:ratio}
 \log\frac{W_{2k}W_rW_1}{W_k^2W_{r+1}}
   =2a(k^2-r)+e,\qquad |e|\leq6c.
\end{equation}
The remarkable feature of the right-hand side is that $k^2-r$ is an integer: unless it is zero, its absolute value is at least one. Hence, if $a>6c$, choose a rational constant $K>1$ such that $6c<\log K<2a-6c$. We call a fixed multiplicative interval of the form $[1/K,K]$ a \emph{height window}\index{height window}. In this case the ratio in \eqref{eq:ratio} belongs to the window exactly when $r=k^2$: for equality the absolute logarithm is at most $6c$, whereas otherwise it is at least $2a-6c$. Writing $K$ as a ratio of positive integers and clearing denominators expresses the test by two polynomial comparisons between products of terms of the sequence.

To illustrate the idea, put $W_j=2^{j^2}$ temporarily (so that there is no error term). For $k=3$ and $r=8,9,10$, the ratio in \eqref{eq:ratio} takes the values $4$, $1$ and $1/4$, respectively: only at $r=9$ does it fall within $[1/2,2]$. A bounded error can be accommodated by spacing the logarithms farther apart and widening the interval. This is merely an illustration: using that sequence itself would presuppose exponentiation. We shall obtain a sequence with property \eqref{eq:model} from the heights of rational points.

\subsection{An elliptic curve as a source of quadratic growth}
Recall that the rational points of an elliptic curve (a nonsingular cubic with distinguished point $\OO$ at infinity\index{point at infinity}) form an Abelian group under the chord--tangent law\index{chord--tangent law},\footnote{For a curve in short Weierstrass form $y^2=x^3+Ax+B$, the law reads as follows. A line meets the cubic in three points, counted with multiplicity, provided $\OO$ is counted as the third intersection of every vertical line. Given $Q_1,Q_2$, let $R$ be the third intersection of the line through $Q_1$ and $Q_2$ (the tangent at $Q_1$ if $Q_1=Q_2$); then $Q_1+Q_2$ is the reflection of $R$ across the $x$-axis. Thus $\OO$ is the neutral element, $-(x,y)=(x,-y)$, and three points add up to $\OO$ exactly when they are collinear. For two affine rational points a nonvertical secant or tangent has rational slope; the third intersection is rational by the cubic equation. Vertical lines and the cases involving $\OO$ are covered by the inverse and identity rules. Associativity is the only property that is not immediate \cite[Chapter III, Section 2]{Silverman}. Explicit formulas are given in Appendix~\ref{app:codes}; Figure~\ref{fig:chord-tangent} illustrates the law on the curve \eqref{eq:curve}.} finitely generated by the Mordell--Weil theorem\index{Mordell--Weil theorem}. To begin with, let us choose a curve whose group of rational points is infinite cyclic:
\begin{equation}\label{eq:curve}
 E:\quad y^2=x^3-74.
\end{equation}
It is a Mordell curve\index{Mordell curve} (one of the form $y^2=x^3+B$ with $B\ne0$), of discriminant $-432\cdot74^2$; its rational-point group has rank one and trivial torsion, and it is generated by\index{generator G@generator $\GG$}
\[
 \GG=(99,985),\qquad 985^2=970\,225=99^3-74.
\]
We shall establish the following group structure:
\begin{equation}\label{eq:group}
 E(\Q)=\{j\GG:j\in\Z\},\qquad j\GG=j'\GG\ \Longleftrightarrow\ j=j'.
\end{equation}
Every rational point of this curve thus has a unique signed index relative to $\GG$. We need not find $\GG$ by an algorithm: its coordinates are fixed constants.

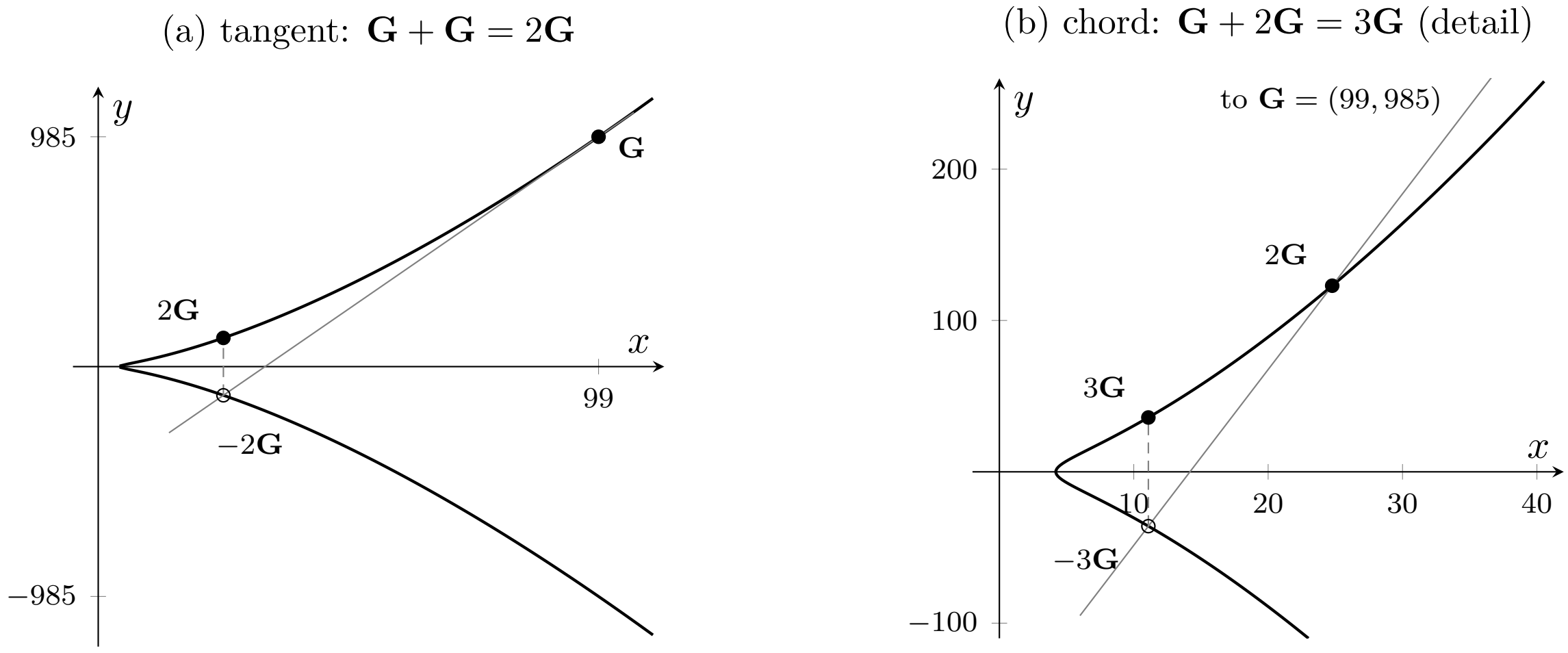
\begin{figure}[ht]
\centering
\begin{tikzpicture}
\begin{axis}[width=.48\linewidth,height=.46\linewidth,axis lines=middle,
  xmin=-5,xmax=112,ymin=-1200,ymax=1200,xtick={99},ytick={-985,985},
  xlabel={$x$},ylabel={$y$},clip mode=individual,tick label style={font=\scriptsize},
  title={(a) tangent: $\GG+\GG=2\GG$},title style={font=\small}]
  \addplot[thick,domain=-1150:1150,samples=300,variable=\t]({(\t*\t+74)^(1/3)},{\t});
  \addplot[gray,domain=14:106,samples=2]{985+14.92538*(x-99)};
  \addplot[densely dashed,gray] coordinates {(24.76699,-122.95595) (24.76699,122.95595)};
  \addplot[only marks,mark=*,mark size=1.8pt] coordinates {(99,985)};
  \addplot[only marks,mark=o,mark size=1.8pt] coordinates {(24.76699,-122.95595)};
  \addplot[only marks,mark=*,mark size=1.8pt] coordinates {(24.76699,122.95595)};
  \node[font=\scriptsize,anchor=west] at (axis cs:100,940) {$\GG$};
  \node[font=\scriptsize,anchor=north] at (axis cs:30,-230) {$-2\GG$};
  \node[font=\scriptsize,anchor=south east] at (axis cs:23,140) {$2\GG$};
\end{axis}
\end{tikzpicture}\hfill
\begin{tikzpicture}
\begin{axis}[width=.48\linewidth,height=.46\linewidth,axis lines=middle,
  xmin=-2,xmax=42,ymin=-110,ymax=260,xtick={10,20,30,40},ytick={-100,100,200},
  xlabel={$x$},ylabel={$y$},clip mode=individual,tick label style={font=\scriptsize},
  title={(b) chord: $\GG+2\GG=3\GG$ (detail)},title style={font=\small}]
  \addplot[thick,domain=-110:258,samples=300,variable=\t]({(\t*\t+74)^(1/3)},{\t});
  \addplot[gray,domain=6:41,samples=2]{985+11.61304*(x-99)};
  \addplot[densely dashed,gray] coordinates {(11.0873,-35.90185) (11.0873,35.90185)};
  \addplot[only marks,mark=*,mark size=1.8pt] coordinates {(24.76699,122.95595)};
  \addplot[only marks,mark=o,mark size=1.8pt] coordinates {(11.0873,-35.90185)};
  \addplot[only marks,mark=*,mark size=1.8pt] coordinates {(11.0873,35.90185)};
  \node[font=\scriptsize,anchor=south east] at (axis cs:24,128) {$2\GG$};
  \node[font=\scriptsize,anchor=north east] at (axis cs:10,-42) {$-3\GG$};
  \node[font=\scriptsize,anchor=south east] at (axis cs:10.5,40) {$3\GG$};
  \node[font=\scriptsize,anchor=east] at (axis cs:34,245) {to $\GG=(99,985)$};
\end{axis}
\end{tikzpicture}
\caption{The chord--tangent law on the real branch of $E:y^2=x^3-74$, which lies over $x\geq74^{1/3}$. (a)~The tangent at $\GG=(99,985)$ meets $E$ again at $-2\GG$; reflection across the $x$-axis gives $2\GG=(96\,118\,209/1970^2,\ 940\,044\,063\,073/1970^3)\approx(24.77,122.96)$. (b)~Detail near the origin: the chord through $\GG$ (outside the frame) and $2\GG$ meets $E$ again at $-3\GG$, whose reflection is $3\GG\approx(11.09,35.90)$. The point $\OO$ cannot be drawn: it is the common third intersection of all vertical lines.}\label{fig:chord-tangent}
\end{figure}

The verification of \eqref{eq:group} for this particular curve is short enough to be given without recourse to a table of elliptic curves. A count of points modulo $7$ and $13$ excludes torsion, and a descent via a rational $3$-isogeny establishes rank one (see Appendix~\ref{app:descent}). We shall see below how a uniform bound on the difference between the logarithmic and canonical heights also establishes that $\GG$, rather than a proper multiple of another point, generates the group.

For a rational number $x=A/B$ in lowest terms with $B>0$, put
\[
 H(x)=\max\{|A|,B\},\qquad h(x)=\log H(x);
\]
these are the \emph{naive} and the \emph{logarithmic height}\index{height!naive}\index{height!logarithmic} of $x$. Both heights measure the size of the reduced numerator and denominator, not the absolute value of the fraction: take, for example, $H(1/1000)=1000$. For an affine point $Q$, write $x(Q)$ and $y(Q)$ for its coordinates and put $h_x(Q)=h(x(Q))$; put $h_x(\OO)=0$.

The canonical (N\'eron--Tate) height\index{height!canonical (N\'eron--Tate)} $\hh$ on an elliptic curve is characterized here by the following familiar facts, with the normalization for the $x$-coordinate \cite[Chapter VIII, Proposition 9.1 and Theorem 9.3]{Silverman}:
\begin{equation}\label{eq:NT}
 \hh(jQ)=j^2\hh(Q),\qquad \hh(Q)>0\text{ if $Q$ is nontorsion},\qquad
 |h_x(Q)-2\hh(Q)|\leq c_E.
\end{equation}
Here $j\in\Z$ and $c_E\geq0$ depends only on the curve $E$. For explicit bounds on height discrepancies in general, see \cite{Silverman1990}. It is the uniformity of the last bound, valid also at $\OO$ with the above convention, that matters for our purposes. For the curve \eqref{eq:curve}, we can estimate the discrepancy rather tightly by examining the duplication formula. For an affine point $Q$, write $x(Q)=A/B$ in lowest terms, $B>0$. For a prime $p$, let $\nu_p$ denote the $p$-adic valuation\index{p-adic valuation@$p$-adic valuation}, that is, the exponent of $p$ in a nonzero rational number (negative when $p$ divides the reduced denominator). Since the equation has integral coefficients, the valuations at every prime dividing $B$ show that $B=d^2$ and $y(Q)=C/d^3$ for an integer $d>0$ and an integer $C$ with $\gcd(A,d)=1$; indeed, a negative valuation of $x$ forces $2\nu_p(y)=3\nu_p(x)$. Thus $C^2=A^3-74d^6$. The duplication formula reads
\[
 x(2Q)=\frac{A(A^3+592B^3)}{4B(A^3-74B^3)}.
\]
Its displayed numerator and denominator are relatively prime. A common prime cannot divide $B$, since $\gcd(A,B)=1$. Apart from $2$, a common prime dividing $A$ must divide $74$; if it does not divide $A$, it must divide $666=592+74$. Hence only $2,3,37$ need be considered. Neither $2$ nor $37$ divides $A$: if $2\mid A$, then $d$ is odd and $C^2\equiv2\pmod4$; if $37\mid A$, then $37\nmid d$ and $\nu_{37}(C^2)=1$. Both conclusions are impossible. As $592$ is divisible by $2$ and $37$, neither of those primes divides the numerator. At $3$, a common divisor would force $A^3\equiv74B^3\pmod3$ and $3\mid C$. If $3\nmid d$, reducing $C^2=A^3-74d^6$ modulo $9$ would then give $A^3\equiv2\pmod9$, which is impossible; if $3\mid d$, then $3\nmid A$ already excludes a common divisor. This proves the asserted coprimality.

All affine real points have $x\geq74^{1/3}>4$. In particular, $H(x(Q))=A$, and the displayed numerator of $x(2Q)$ is between $A^4$ and $9A^4$, whereas its positive denominator is less than $4BA^3<A^4$. Consequently
\[
 0\leq h_x(2Q)-4h_x(Q)\leq\log9.
\]
The same bounds hold at $\OO$ by definition. Summation along $Q,2Q,4Q,\ldots$, and the canonical-height limit $2\hh(Q)=\lim_{s\to\infty}4^{-s}h_x(2^sQ)$, give
\begin{equation}\label{eq:one-sided-height}
 0\leq2\hh(Q)-h_x(Q)\leq\frac{\log9}{3}=c_E,\qquad c_E=\frac23\log3.
\end{equation}
Here the limit and quadraticity of $\hh$ are the classical canonical-height theorem; the explicit discrepancy follows from the elementary estimates just given.

We can now finish the verification of \eqref{eq:group}. If $\GG$ were a multiple $jQ$ of a generator of the infinite cyclic group, with $|j|\geq2$, then \eqref{eq:one-sided-height} and $x(Q)>4$ would imply $2\hh(\GG)=j^2\cdot 2\hh(Q)\geq4\log5$. Yet $2\hh(\GG)\leq\log99+\frac23\log3<4\log5$ (indeed $99\cdot3^{2/3}<297<625$). Hence $|j|=1$, as desired.

Define
\begin{equation}\label{eq:Wdef}
 W(Q)=\begin{cases}H(x(Q)),&Q\ne \OO,\\1,&Q=\OO.\end{cases}
 \qquad W_j=W(j\GG)\quad(j\geq0).
\end{equation}
Since $x(-Q)=x(Q)$, $W(j\GG)=W_{|j|}$. By \eqref{eq:one-sided-height},
\begin{equation}\label{eq:W}
 W_0=1,\qquad \log W_j=aj^2+\varepsilon_j,\qquad
 a=2\hh(\GG),\quad -c_E\leq\varepsilon_j\leq0.
\end{equation}
At $j=0$ take $\varepsilon_0=0$. The one-sided bound improves the estimate in \eqref{eq:ratio}: the three positive and the three negative occurrences of the errors (counting repetitions) give $|e|\leq3c_E$. We may therefore choose the fixed integer $K=10$ since
\begin{equation}\label{eq:constants}
 a>c_E,\qquad 3c_E<\log K<2a-3c_E.
\end{equation}
Indeed $a\geq\log99$, $3c_E=\log9$, and $2\log99-\log9=\log1089>\log10$. No amplification of the height\index{height amplification} by a fixed multiple of the point is needed for this curve; $W(\GG)=99$. The sequence $W_j$ is strictly increasing since the difference of consecutive logarithms is at least $a(2j+1)-c_E>0$, and it is unbounded.

Henceforth we shall apply arithmetic operations directly to rational points and compare their heights. Each point is represented by a unique natural tuple encoding its reduced rational coordinates and distinguishing the point at infinity. Appendix~\ref{app:codes} shows how to do so by univocal polynomial specifications on natural codes. Only group additions and comparisons of natural heights are needed in the construction below; signed indices are used in its justification, not in its specification.

\subsection{Orienting the indices without extra sign choices}
At first sight the symmetry $W(j\GG)=W_{|j|}$ is an obstacle to the univocal selection of a point: its height cannot discriminate between two opposite indices. By comparing the heights of two consecutive points, however, we can single out the nonnegative indices:
\begin{equation}\label{eq:orient}
 W(Q+\GG)>W(Q)\quad\Longleftrightarrow\quad Q=j\GG\text{ for some }j\geq0.
\end{equation}
Indeed, if $Q=j\GG$ with $j\geq0$, the absolute index increases by one on adding $\GG$; if $j<0$, it decreases by one. The points $Q$ satisfying \eqref{eq:orient} constitute a set which we denote by $\ord$\index{oriented set}. It includes the identity, and for $Q=j\GG\in\ord$ we have $W(Q+\GG)=W_{j+1}$. This comparison replaces the selection among torsion translates and the maxima of two neighbouring heights used in \cite{OpenAI2026}.

\section{A singlefold function of intermediate growth}\label{sec:growth}
We are now in a position to realize the idea suggested by \eqref{eq:ratio}: given a number $u$, locate it between successive heights and replace the index thus singled out by its square. The height at the latter index will provide the growth function needed in the Pell construction.
\begin{proposition}\label{prop:growth}\index{growth function}
There is a total function $g:\N\to\N_{>0}$ with a singlefold polynomial graph such that
\begin{equation}\label{eq:asymp}
 \log g(u)\sim\frac{(\log u)^2}{a}\qquad (u\longrightarrow\infty),
\end{equation}
where $a>0$ is a fixed constant. In particular, $g$ grows faster than every fixed power of $u$ and more slowly than $u^u$.
\end{proposition}
\begin{proof}
Given $u\in\N$, choose the rational point $Q$ by the bracket\index{bracket}
\begin{equation}\label{eq:bracket}
 W(Q)\leq u+1<W(Q+\GG).
\end{equation}
Taken together, the two inequalities imply \eqref{eq:orient}, whence $Q=k\GG$ for some $k\geq0$. There is exactly one such $k$, since the increasing sequence of positive integers $W_j$ starts at 1 and is unbounded: it is the one for which $W_k\leq u+1<W_{k+1}$. Thus \eqref{eq:bracket} determines $Q$ uniquely.

The next step is to select in $\ord$ the point whose index is the square of the index of $Q$. To this end, require $R\in\ord$ by \eqref{eq:orient}, so that $R=r\GG$ with $r\geq0$ (the index $r$ is used solely in the proof), and form the two positive integers
\[
 L=W(2Q)W(R)W(\GG),\qquad D=W(Q)^2W(R+\GG).
\]
Impose the two inequalities
\begin{equation}\label{eq:window}
 L\leq KD,\qquad D\leq KL.
\end{equation}
The ratio $L/D$ belongs to the height window $[1/K,K]$ just when $|\log(L/D)|\leq\log K$. Now \eqref{eq:ratio} applies also when either index is zero, since $\varepsilon_0=0$; thus
\[
 \log(L/D)=2a(k^2-r)+e,\qquad |e|\leq3c_E.
\]
When $r=k^2$, the error term alone contributes to the logarithm, and $3c_E<\log K$ ensures that \eqref{eq:window} holds. For any other integral $r$ the absolute value of the logarithm is at least $2a-3c_E>\log K$, contrary to one of the two inequalities. Hence the sole choice of $R$ is $k^2\GG$.

Set $g(u)=W(R)=W_{k^2}$. The two unknown rational points have now been determined uniquely, and each arithmetic call contributing to \eqref{eq:bracket} and \eqref{eq:window} has a unique complete internal tuple by Appendix~\ref{app:codes} and Lemma~\ref{lem:compiler}. These conditions are finite polynomial constraints on their natural codes: the only multiplications of points are by fixed integers, every height is an integer, and the displayed comparisons have unique slacks. The proof indices $k,r$ are not quantified variables. Thus the graph of $g$ is specified univocally. This reasoning includes $u=0$, where $k=r=0$ and $g(0)=1$.

As $u\to\infty$, the index $k$ tends to infinity, and \eqref{eq:bracket} gives
\[
 ak^2-c_E\leq\log(u+1)<a(k+1)^2+c_E.
\]
Dividing by $ak^2$ and using $(k+1)^2/k^2\to1$, we get $\log(u+1)\sim ak^2$. Inserting this in $\log g(u)=ak^4+O(1)$ proves \eqref{eq:asymp}. The two growth comparisons follow by dividing $\log g(u)$ by $\log u$ and by $u\log u$, respectively.
\end{proof}

For the Pell construction we will use the following consequence of \eqref{eq:asymp}: there is a fixed $\kappa>0$ such that for all sufficiently large $u$
\begin{equation}\label{eq:g-bounds}
 \kappa(\log u)^2<\log g(u)<(u-1)\log u.
\end{equation}
The lower bound will leave room for the desired Pell solutions, whereas the upper bound will rule out solutions whose indices are too large. Both are estimates used to justify a polynomial specification, not constraints involving logarithms to be adjoined to it.

For the curve \eqref{eq:curve} the constants can be made explicit: \eqref{eq:g-bounds} holds with
\begin{equation}\label{eq:explicit-delta}
 \kappa=\frac1{32\log206}\approx5.87\cdot10^{-3}\qquad\text{for every }u\geq10^5.
\end{equation}
Indeed, by \eqref{eq:W} with $j=1$ we have $\log99\leq a\leq\log99+c_E=\log(99\cdot3^{2/3})<\log206$. Let $u\geq10^5$ and $L=\log(u+1)$. Then $u+1\geq W_1$, so the index $k$ singled out by \eqref{eq:bracket} is at least $1$, and $L<\log W_{k+1}\leq a(k+1)^2\leq4ak^2$. Hence $\log g(u)\geq ak^4-c_E>L^2/(16a)-c_E\geq2\kappa(\log u)^2-c_E$, which exceeds $\kappa(\log u)^2$ because $(\log u)^2\geq(\log10^5)^2>132>c_E/\kappa$. On the other hand $ak^2\leq L+c_E$, so that $\log g(u)\leq ak^4\leq(L+c_E)^2/\log99$. Since $L+c_E=\log(u+1)+c_E<2\log u$ for $u\geq10^5$ (use $c_E=\frac23\log3$), we obtain $\log g(u)<4(\log u)^2/\log99<(u-1)\log u$, using $\log99>4$ and $\log u\leq u-1$.

\section{Recovering exponentiation from Pell equations}\label{sec:pell}
\subsection{Exponentiation as a quotient of Pell solutions}
Let us first see why the search for a polynomial specification of exponentiation brings Pell equations into play. For an integer $A\geq2$, let $x_i(A),y_i(A)$ be defined by
\begin{equation}\label{eq:pell-sequence}
 x_i(A)+y_i(A)\sqrt{A^2-1}=(A+\sqrt{A^2-1})^i\quad(i\geq0).
\end{equation}
They are integer solutions of $x^2-(A^2-1)y^2=1$. Their second coordinates satisfy
\begin{equation}\label{eq:recurrence}
 y_0(A)=0,\quad y_1(A)=1,\quad y_{i+2}(A)=2Ay_{i+1}(A)-y_i(A).
\end{equation}
The first few terms are
\[
 \begin{array}{c|l}
 i&y_i(A)\\\hline
 0&0\\
 1&1\\
 2&2A\\
 3&4A^2-1\\
 4&8A^3-4A\\
 5&16A^4-12A^2+1.
 \end{array}
\]
For fixed $n\geq0$, the recurrence shows that $y_{n+1}(A)$ is a polynomial in $A$ of degree $n$ with leading coefficient $2^n$. It follows that
\begin{equation}\label{eq:pell-limit}
 \lim_{u\to\infty}\frac{y_{n+1}(2u)}{y_{n+1}(u)}=2^n.
\end{equation}
For example, with $n=3$ the ratio is exactly $8+6/(2u^2-1)$:
\begin{center}
\begin{tabular}{rrrr}
\toprule
$u$&$y_4(u)$&$y_4(2u)$&Ratio (rounded)\\
\midrule
2 & 56 &496&8.857143\\
4 &496&4064&8.193548\\
8 &4064&32704&8.047244\\
16&32704&262016&8.011742\\
\bottomrule
\end{tabular}
\end{center}
Observe that the integer part is $8$ in all four instances. This observation underlies the use of paired Pell equations to specify exponentiation (cf. \cite{Robinson1952,Mat1974,CCO2024}).

Three difficulties remain before we can make such a specification univocal. We must single out the $(n+1)$st solutions of the two Pell equations without using an index operator; we must determine the scale $u$ rather than simply require it to be large enough; and we must ensure that the quotient lies between $2^n$ and $2^n+1$, the latter excluded. The new growth function will attend to the first issue, while a constraint involving the proposed value of the power will attend to the other two.

\subsection{Classical facts about Pell sequences}
\begin{lemma}\label{lem:pell}\index{Pell equation}
For each integer $A\geq2$, all solutions in positive integers of $x^2-(A^2-1)y^2=1$ are uniquely of the form $(x_i(A),y_i(A))$ for $i\geq1$. For $i\geq1$ one has
\begin{align}
 y_i(A)&\equiv i\pmod{A-1},\label{eq:congruence}\\
 (2A-1)^{i-1}&\leq y_i(A)\leq(2A)^{i-1}.\label{eq:pell-bounds}
\end{align}
\end{lemma}
\begin{proof}
Put $D=A^2-1$. It lies strictly between $(A-1)^2$ and $A^2$, so $\sqrt D$ is irrational. The fundamental positive solution is $(A,1)$. For completeness, let $\alpha=A+\sqrt D$ and let $\beta=x+y\sqrt D$ come from any positive solution. Choose $i\geq0$ so that $1\leq\gamma=\beta\alpha^{-i}<\alpha$. Since $\alpha^{-1}=A-\sqrt D$, the number $\gamma$ belongs to $\Z[\sqrt D]$ and has norm one. Write $\gamma=x'+y'\sqrt D$. If $\gamma>1$, then $y'=(\gamma-\gamma^{-1})/(2\sqrt D)>0$ and $x'=(\gamma+\gamma^{-1})/2>0$. Hence $y'\geq1$ and $x'^2=1+Dy'^2\geq A^2$, giving $\gamma\geq\alpha$, impossible. Thus $\gamma=1$, $\beta=\alpha^i$, and $i\geq1$ is unique. This proves the first assertion, and \eqref{eq:recurrence} follows from the identity $\alpha+\alpha^{-1}=2A$.

Modulo $A-1$, the number $A$ is congruent to 1, and the recurrence gives $y_{i+2}\equiv2y_{i+1}-y_i$; the initial values prove \eqref{eq:congruence} by induction. To obtain the growth estimates, let $R_i(A)=y_{i+1}(A)/y_i(A)$ for $i\geq1$. Then
\[
 R_1(A)=2A,\qquad R_{i+1}(A)=2A-\frac1{R_i(A)}.
\]
Induction yields $2A-1\leq R_i(A)\leq2A$. Multiplying these bounds and using $y_1(A)=1$ proves \eqref{eq:pell-bounds}.

\end{proof}

\subsection{When the quotient is exact}
To turn the suggestive limit \eqref{eq:pell-limit} into a valid biimplication, we shift the argument of the numerator and bound the scale in terms of the proposed output. Write $A\mathbin{\mathrm{div}}B=\lfloor A/B\rfloor$ when $A$ is nonnegative and $B$ a positive integer.

\begin{lemma}[Matiyasevich's quotient formula]\label{lem:quotient}\index{Matiyasevich's quotient formula}
Let $b,n,c,u\in\N$, with $b\geq1$ and $u\geq3(c+1)(n+1)$. Then
\begin{equation}\label{eq:quotient}
 b^n=c\quad\Longleftrightarrow\quad
 c=y_{n+1}(bu+1)\mathbin{\mathrm{div}}y_{n+1}(u).
\end{equation}
\end{lemma}
This is the second assertion of \cite[Lemma 6.1, p.~600]{CCO2024}, where it is traced to \cite[p.~308 of the English translation]{Mat1974}. Its usefulness for the present argument is that the same bound, depending on the proposed output, works in both directions. Although our immediate concern is with base 2, we shall exploit the formula for arbitrary positive bases when several power equations share a cutoff. Appendix~\ref{app:quotient} gives an elementary proof of both implications.

For example, take $b=2$, $n=3$ and $c=8$. The prescribed lower bound for $u$ is $108$. At $u=108$ we get
\[
 y_4(108)=10\,077\,264,\qquad y_4(217)=81\,745\,636,
\]
and the integer quotient is indeed $8$. The shift $+1$ in the numerator is part of the classical formula; it was absent from the illustrative limit \eqref{eq:pell-limit} and from the ratio used in \cite[Lemma 5.1 and Proposition 5.2]{OpenAI2026}.

The salient point about \eqref{eq:quotient} is that we need not know in advance whether $c=b^n$: a wrong proposed value of $c$ cannot satisfy the quotient condition. What remains to be shown is how a suitable upper bound on the Pell coordinates permits us to select precisely the solutions entering that quotient.

\subsection{The paired-Pell system with a prescribed cutoff}
The index-selection and unique-completion argument in \cite[Proposition 5.2]{OpenAI2026} can be separated from the choice of its growth function. In the following version, based on the classical shifted quotient, both the scale and the cutoff are supplied, and the exact admissibility condition is recorded. We call $u$ the \emph{scale}\index{scale}, since it determines the Pell parameters, and $v$ the \emph{cutoff}\index{cutoff}\index{cutoff!admissible}, a strict upper bound on the second coordinates. Within the hypotheses of the following lemma, a cutoff is admissible for a true power equation if it exceeds both prescribed second coordinates.
\begin{lemma}[Exact cutoff lemma]\label{lem:exact-cutoff}\index{exact cutoff lemma}
Fix $b,n,c,u,v\in\N$ with $b\geq1$, $u\geq3(c+1)(n+1)$ and
\[
 1\leq v<(2u-1)^{u-1}.
\]
Consider the system in the seven natural unknowns $X,Y,F,H,q,q',r$:
\begin{equation}\label{eq:general-pell}
\begin{gathered}
 X^2-((bu+1)^2-1)F^2=1,\quad Y^2-(u^2-1)H^2=1,\\
 F<v,\quad H<v,\\
 F=n+1+buq,\quad H=n+1+(u-1)q',\\
 F=cH+r,\quad r<H.
\end{gathered}
\end{equation}
It has exactly one solution, including all comparison slacks, if and only if
\begin{equation}\label{eq:admissible}
 c=b^n,\qquad v>\max\{y_{n+1}(bu+1),y_{n+1}(u)\}.
\end{equation}
Otherwise it has no solution.
\end{lemma}
\begin{proof}
Let $(X,F)$ and $(Y,H)$ solve the two Pell equations of \eqref{eq:general-pell}. The congruence equations in \eqref{eq:general-pell} already guarantee that $F,H\geq n+1>0$; no separate positivity constraints are needed. By Lemma~\ref{lem:pell}, their positive second coordinates have indices $i,j\geq1$. Were either index at least $u$, the corresponding coordinate would be at least $(2u-1)^{u-1}$, by \eqref{eq:pell-bounds} and $bu+1\geq u$, in contradiction to the upper cutoff. The congruences give $i\equiv n+1\pmod{bu}$ and $j\equiv n+1\pmod{u-1}$. Since $1\leq i,j\leq u-1$ and $1\leq n+1<u-1$, both indices equal $n+1$ (for $b=1$ the first modulus is $u$, which is still sufficient). The quotient formula in Lemma~\ref{lem:quotient} forces $c=b^n$, and the inequalities give the second requirement in \eqref{eq:admissible}. The endpoint $j=u-1$ is inconsequential here: its residue is zero rather than $n+1$.

Conversely, under \eqref{eq:admissible} choose the prescribed Pell coordinates and their positive first coordinates. Their congruences give uniquely determined integral $q,q'$, nonnegative because $y_{n+1}(A)\geq(2A-1)^n\geq3^n\geq n+1$. The first coordinates $X,Y$ are the unique positive natural roots of the respective Pell equations. The quotient formula gives the unique remainder $r$. These values satisfy the whole system. In either direction all the coordinates and the three comparison slacks are determined once the displayed arguments are fixed. The numerical upper bound involving a power belongs to the hypotheses, not to the polynomial specification \eqref{eq:general-pell}.
\end{proof}

\subsection{The polynomial system for powers of two}
\begin{proposition}[Powers of two]\label{prop:pow}\index{powers of two}
Suppose a total function $g:\N\to\N_{>0}$ has a singlefold polynomial graph and satisfies \eqref{eq:g-bounds} for some $\kappa>0$ and all $u\geq u_0$, where $u_0\geq5$ is a fixed integer. Then
\[
 \Pow(n,c)\quad\Longleftrightarrow\quad c=2^n
\]
has a singlefold polynomial representation.
\end{proposition}
\begin{proof}
Fix an integer $\ell>2/(\kappa\log2)$; for the function $g$ of Proposition~\ref{prop:growth}, with the values \eqref{eq:explicit-delta}, one may take $u_0=10^5$ and $\ell=492$, since $2/(\kappa\log2)=64\log_2206<491.95$. Consider the following system of conditions in natural variables, with $n,c$ the displayed arguments:
\begin{align}
 u&=3(c+1)^\ell(n+1)+u_0,& v&=g(u),\label{eq:scale}\\
 X^2-((2u+1)^2-1)F^2&=1,&Y^2-(u^2-1)H^2&=1,\label{eq:paired}\\
 F&<v,&H&<v,\label{eq:cutoff}\\
 F&=n+1+2u q,&H&=n+1+(u-1)q',\label{eq:residues}\\
 F&=cH+r,&r&<H.\label{eq:division}
\end{align}
This is polynomial once the singlefold system for $g$ has been substituted and comparisons have been expanded. In particular $\ell$ is a \emph{fixed} integer exponent, not an exponent variable.

\smallskip\noindent\emph{Soundness.} For any candidate $c$, the polynomial scale equation ensures $u\geq3(c+1)(n+1)$ and $u\geq u_0$. By \eqref{eq:g-bounds}, $1\leq v=g(u)<u^{u-1}<(2u-1)^{u-1}$. Thus Lemma~\ref{lem:exact-cutoff}, with $b=2$, applies to the remaining equations \eqref{eq:paired}--\eqref{eq:division}. It proves immediately that no false candidate $c$ has a witness, and that a correct candidate can have at most one complete Pell tuple. We see here why the scale equation involves $c$: the bound required by the quotient formula holds before we assume $c=2^n$.

\smallskip\noindent\emph{Completeness.} Suppose $c=2^n$ and take the unique $u,v$ and complete witness for the call $v=g(u)$. The outstanding issue is whether the two prescribed Pell coordinates actually fall below the cutoff $v$. From the scale we have $\log u\geq\ell n\log2$. Since $u\geq5$, $4u+2\leq u^2$. By \eqref{eq:pell-bounds},
\[
 \log y_{n+1}(2u+1),\ \log y_{n+1}(u)
 \leq n\log(4u+2)
 \leq\frac{2(\log u)^2}{\ell\log2}
 <\kappa(\log u)^2<\log g(u).
\]
Consequently the second condition of \eqref{eq:admissible} holds. The exact cutoff lemma supplies the remaining tuple.

\smallskip\noindent\emph{Uniqueness.} The scale equation fixes $u$, the singlefold graph of $g$ fixes $v$ and its internal witnesses, and Lemma~\ref{lem:exact-cutoff} fixes all seven Pell variables and their slacks. This includes $n=0$, when the Pell coordinates are $F=H=1$, $X=2u+1$, $Y=u$ and $q=q'=r=0$.
\end{proof}

\section{Completion of the proof}\label{sec:conclusion}
Proposition~\ref{prop:growth} supplies the function $g$ in Proposition~\ref{prop:pow}. The latter gives a singlefold polynomial representation of the graph of $2^u=v$. There is thus a polynomial $F$ with integer coefficients such that $2^u=v$ holds if and only if there is exactly one tuple $\vz\in\N^s$ satisfying $F(u,v,\vz)=0$, and no such tuple otherwise. The passage from a finite system to $F$ is the sum-of-squares construction in Section~\ref{sec:systems}.

For an arbitrary recursively enumerable set $S\subseteq\N^n$, take the polynomial $D$ provided by Theorem~\ref{thm:mat}. The polynomial
\[
 P(\va,\mathbf{x},u,v,\vz)=D(\va,\mathbf{x},u,v)^2+F(u,v,\vz)^2
\]
is the desired representation by the substitution principle following Theorem~\ref{thm:mat}. This proves Theorem~\ref{thm:main}. The finitefold assertion follows immediately.\hfill$\square$

For the fixed curve the construction uses the generator $\GG=(99,985)$, the window constant $K=10$, no height amplification, the growth constant $\kappa=1/(32\log206)$, the scale threshold $u_0=10^5$ and the exponent $\ell=492$. The growth constant is used in the estimates; the actual polynomial system uses only fixed integer constants and natural codes. These data specify the arithmetic module, not the coefficients of the entire representing polynomial for an arbitrary r.e. set: the latter also uses the polynomial $D$ from Theorem~\ref{thm:mat}.

\part{General principles and further results}\label{part:extensions}
The constructions to be developed next stem from two different ingredients of Part~\ref{part:proof}. The height comparisons will enable us to extend the arithmetic of indices and to pass to a cyclic quotient even in the presence of higher rank. The selection of Pell solutions will, in turn, let us relate the multiplicity of a specification of exponentiation to that of a supplied growth relation. The two lines of argument can be read independently; their essential dependencies are indicated in the following table.
\begin{center}
\begin{tabular}{@{}>{\raggedright\arraybackslash}p{.46\linewidth}>{\raggedright\arraybackslash}p{.46\linewidth}@{}}
\toprule
Height constructions & Cutoff constructions\\\midrule
Height estimate \eqref{eq:W} and window \eqref{eq:window} $\longrightarrow$ index arithmetic and hierarchy & Lemma~\ref{lem:exact-cutoff} $\longrightarrow$ Theorems~\ref{thm:transfer} and \ref{thm:joint-transfer}\\[2pt]
Theorem~\ref{thm:rank-quotient} $\longrightarrow$ Theorem~\ref{thm:abstract-height} $\longrightarrow$ cutoff & Transfer $\longrightarrow$ sparse-domain, one-call and iterated criteria\\[2pt]
Quadratic-height cutoff plus Proposition~\ref{prop:pow} gives powers of two. & Theorem~\ref{thm:mat} passes the witness bound to r.e. relations.\\
\bottomrule
\end{tabular}
\end{center}
Except where explicitly stated otherwise, the constructions in this part use only the supplied polynomial systems and the operations specified in their proofs, without invoking Theorem~\ref{thm:main}. We call these constructions \emph{relative}\index{relative construction}, without introducing a formal notion of reducibility.

\section{Arithmetic of indices and an abstract height principle}\label{sec:indices-principle}
\subsection{Addition and multiplication of elliptic indices}
Selecting the square of an index, as we did in Proposition~\ref{prop:growth}, suggests a more general use of the same height comparisons. Relative to the oriented set $\ord$ of \eqref{eq:orient}, the bijection $i\mapsto i\GG$ from $\N$ onto $\ord$ already identifies addition with the group law. The polarization of the quadratic term in a height estimate will enable us to specify multiplication as well.
\begin{proposition}\label{prop:indices}
On the fixed canonical natural codes of $\ord$, the relations
\[
 \operatorname{Add}(i\GG,j\GG,k\GG)\iff i+j=k,
 \qquad \operatorname{Mult}(i\GG,j\GG,k\GG)\iff ij=k
\]
have singlefold polynomial systems. Their arguments and outputs are point codes, not the usual natural indices.
\end{proposition}
\begin{proof}
Check that the input and output points belong to $\ord$ via \eqref{eq:orient}; these tests have unique witnesses for fixed points. Addition is the group law. For multiplication, write $Q=i\GG$, $R=j\GG$, $S=k\GG$ with $i,j,k\geq0$, and put
\[
 L=W(Q+R)W(S)W(\GG),\qquad D=W(Q)W(R)W(S+\GG).
\]
By \eqref{eq:W},
\begin{align*}
 \log(L/D)&=a\big((i+j)^2+k^2+1-i^2-j^2-(k+1)^2\big)+e\\
 &=2a(ij-k)+e,\qquad |e|\leq3c_E.
\end{align*}
The window $L\leq KD$, $D\leq KL$ therefore holds exactly when $ij=k$. Every arithmetic call has unique witnesses, as before. Because $\ord$ contains one point of each nonnegative index, the output is uniquely determined when both input points are fixed.
\end{proof}
This construction extends to any rank-one elliptic curve over $\Q$ with trivial rational torsion (after putting it in short Weierstrass form). With a two-sided bound on the height discrepancy, one uses a fixed scalar multiple and a larger window as in Theorem~\ref{thm:abstract-height} below; the one-sided estimate for our Mordell curve is not needed for that extension. For curves with torsion one may use the canonical quotient construction of Corollary~\ref{cor:torsion-height}. The proposition gives an interpretation of arithmetic using \emph{natural codes of rational points}, including reduced numerators and denominators. It does not supply a Diophantine interpretation over the field $\Q$ itself; in particular, it is not a solution of Hilbert's tenth problem over $\Q$.

\begin{corollary}[A hierarchy of height functions]\label{cor:hierarchy}\index{hierarchy of height functions}
Put $a=2\hh(\GG)$. For every fixed integer $d\geq1$, there is a total function $g_d:\N\to\N_{>0}$ with a singlefold polynomial graph such that
\[
 \log g_d(u)\sim\frac{(\log u)^d}{a^{d-1}}\qquad (u\to\infty).
\]
For $d\geq2$, these functions grow faster than every fixed power of $u$ and more slowly than $u^u$. For integers $1\leq d<d'$, $\log g_d(u)/\log g_{d'}(u)\to0$.
\end{corollary}
\begin{proof}
Select $Q=k\GG$ by \eqref{eq:bracket}. Starting with $Q_1=Q$, use the multiplication system of Proposition~\ref{prop:indices} to require $\operatorname{Mult}(Q_j,Q,Q_{j+1})$ for $1\leq j<d$. Its output is forced to be $Q_{j+1}=k^{j+1}\GG$, including when $k=0$. The number of calls is fixed since $d$ is fixed; for $d=1$ there are none. Set $g_d(u)=W(Q_d)=W_{k^d}$. Each intermediate point and every auxiliary tuple are unique by the preceding proposition and Lemma~\ref{lem:compiler}. The bracket gives $\log u\sim ak^2$, whereas \eqref{eq:W} gives $\log g_d(u)=ak^{2d}+O(1)$. Combining the two yields the displayed equivalent. The remaining comparisons follow by division by $\log u$, by $u\log u$, or by $\log g_{d'}(u)$, respectively.
\end{proof}
The original growth function is $g_2$. For $d=1$ the logarithmic equivalent alone gives no claim of superpolynomial growth. The corollary concerns fixed exponents $d$, not a single graph with $d$ as an additional input.

\subsection{An abstract reduction to quadratic logarithmic heights}\label{sec:height-principle}
In this abstract setting a height is simply a positive integer-valued function $H$ on a naturally coded group. Along the multiples $jP$ ($j\in\Z$) of a fixed group element $P$, we call it a polynomial logarithmic height\index{height!polynomial logarithmic} if $\log H(jP)$ differs from a real polynomial in $j$ by a uniformly bounded amount. Once we have the estimate \eqref{eq:W}, very little of the geometry of the curve enters the subsequent constructions. Even the quadratic main term need not be given at the outset: a fixed finite difference reduces a polynomial main term of degree at least two to degree two. Observe that only even degrees can occur under a bounded-error estimate on all signed indices, since $H$ is positive integer-valued whereas a polynomial of odd degree with positive leading coefficient tends to minus infinity when $j\to-\infty$.

\begin{theorem}[Polynomial logarithmic-height reduction principle]\label{thm:abstract-height}\index{polynomial logarithmic-height reduction principle}
Let $G$ be an infinite cyclic group with identity $O$ and a fixed generator $P$. Suppose that its elements have unique codes in $\N^s$, that validity of a code, addition and negation have singlefold polynomial systems on these codes, and that $H:G\to\N_{>0}$ has a singlefold polynomial graph. Assume that there are a real polynomial $p$ of even degree $d\geq2$ with positive leading coefficient $\alpha$, and a constant $C\geq0$, such that
\begin{equation}\label{eq:polynomial-height}
 \bigl|\log H(jP)-p(j)\bigr|\leq C\qquad(j\in\Z).
\end{equation}
Then the following systems can be constructed using only the assumed systems, finitely many fixed group operations, and elementary polynomial constraints:
\begin{enumerate}
\item a positive integer-valued height $W$ with a singlefold polynomial graph, satisfying $W(O)=1$, $W(-Q)=W(Q)$ for every $Q\in G$, and $\log W(jP)=aj^2+O(1)$, with the quadratic coefficient $a$ arbitrarily large while the error bound stays fixed;
\item a singlefold polynomial specification of the domain of codes corresponding to the nonnegative indices and singlefold systems for addition and multiplication of these indices;
\item a total singlefold function $g:\N\to\N_{>0}$ with $\log g(u)\sim (\log u)^2/a$.
\end{enumerate}
In particular the hypotheses suffice for the singlefold representation theorem via Proposition~\ref{prop:pow}. The construction in (1)--(3) does not use exponentiation with a variable exponent or that representation theorem.
\end{theorem}
\begin{proof}
Put $h=d-2$, and define a positive rational-valued\index{finite differences} function by the fixed product
\[
 D(Q)=\prod_{i=0}^{h}H(Q+iP)^{(-1)^{h-i}\binom hi}.
\]
A negative exponent means a factor in the denominator, and every exponent here is a fixed integer. Thus $D$ can be computed as a numerator and positive denominator, without needing a new rational primitive. For $h=0$ the product is just $H(Q)$. If $\Delta f(j)=f(j+1)-f(j)$, then
\[
 \log D(jP)=\Delta^h p(j)+e_j,\qquad |e_j|\leq 2^h C.
\]
The polynomial $q(j)=\Delta^h p(j)$ has degree two and leading coefficient $A=\alpha d!/2$. Symmetrization and normalization remove its linear and constant terms exactly:
\begin{align*}
 S(Q)&=D(Q)D(-Q)/D(O)^2,\\
 \log S(jP)&=q(j)+q(-j)-2q(0)+e_j+e_{-j}-2e_0\\
 &=\alpha d!\,j^2+\eta_j,\qquad |\eta_j|\leq 2^d C.
\end{align*}
In particular $S(O)=1$ and $S(-Q)=S(Q)$. Choose a fixed positive integer $m$\index{height amplification} with $\alpha d!m^2>2^d C$. Then $S(mQ)\geq1$ for every $Q$ (strictly greater than one for $Q\ne O$). Define
\[
 W(Q)=\lceil S(mQ)\rceil,\qquad a=\alpha d!m^2.
\]
Since $1\leq S(mQ)\leq W(Q)<S(mQ)+1\leq2S(mQ)$, the logarithmic rounding error is at most $\log2$. Hence
\begin{equation}\label{eq:abstract-error}
 W(O)=1,\quad W(-Q)=W(Q),\quad
 |\log W(jP)-aj^2|\leq C_*:=2^d C+\log2.
\end{equation}
The bound $C_*$ does not depend on $m$. The ceiling is itself singlefold: for a computed fraction $N/B$ with $N,B>0$, require $N=B(w-1)+r$ with $w\geq1$ and $1\leq r\leq B$. This fixes $w,r$ and their comparison slacks uniquely. The fraction need not be reduced because its numerator and denominator are computed products, not free choices. All the other operations are fixed calls to the given systems or fixed products. Thus $W$ has a singlefold polynomial graph.

Choose an integer $K>\exp(6C_*)$ and increase $m$ if needed so that $a>2C_*$ and $2a-6C_*>\log K$. Now $W_j=W(jP)$, $j\geq0$, is strictly increasing and unbounded by the same difference estimate as in Section~\ref{sec:height}. The inequality $W(Q+P)>W(Q)$ selects precisely the nonnegative indices. On that domain addition is the group law, and for $Q=iP$, $R=jP$, $T=kP$ the ratio
\[
 \frac{W(Q+R)W(T)W(P)}{W(Q)W(R)W(T+P)}
\]
has logarithm $2a(ij-k)$ with error at most $6C_*$. The two comparisons for the fixed window $[1/K,K]$ therefore select exactly the product index $k=ij$, with a unique complete tuple. This proves (2), including zero indices. Finally choose $Q=kP$ by $W(Q)\leq u+1<W(Q+P)$; use the multiplication system to select $T=k^2P$ and set $g(u)=W(T)$. The bracket and the estimate \eqref{eq:abstract-error} give $\log u\sim ak^2$ and $\log g(u)\sim ak^4$, exactly as in Proposition~\ref{prop:growth}. This proves (3). No part of the proof uses a uniform operation $(j,Q)\mapsto jQ$: the sole scalar $m$ and all shifts are fixed.
\end{proof}

For the integer $x$-height of a rank-one elliptic curve the polynomial main term of its logarithm already has degree two, so the finite differences are not needed; the direct construction in Section~\ref{sec:height} is simpler. The theorem separates the algebraic source of a height from the arithmetic consequences of its growth. One must still supply a genuine coded group and a height satisfying \eqref{eq:polynomial-height}; an asymptotic equivalent without bounded additive error is not sufficient for this proof. The theorem does not itself treat higher-rank groups; the next constructions address this issue. The reduction is structural: the customary canonical heights on Abelian varieties associated with symmetric ample line bundles are quadratic; finding a natural higher-degree source remains open (see Section~\ref{sec:cco-programme}).

\subsection{Finite torsion and the quotient by a finite subgroup}
The passage from an infinite cyclic group to one which is cyclic modulo finite torsion calls for a canonical choice of representative in each torsion coset. Since here each coset is a finite list of points, the difficulty is essentially one of univocal coding, independent of the particular curve.
\begin{corollary}[Virtually cyclic height principle]\label{cor:torsion-height}\index{virtually cyclic height principle}
Let $G=\Z P\oplus T$, where $T$ is a fixed finite Abelian group. Assume unique natural codes for $G$, singlefold validity, addition, negation and a singlefold positive integer-valued height $H$, as in Theorem~\ref{thm:abstract-height}. Suppose
\[
 |\log H(jP+t)-p(j)|\leq C\qquad(j\in\Z,\ t\in T)
\]
for some constant $C\geq0$ and some real polynomial $p$ of even degree at least two with positive leading coefficient. Then $G/T$ has canonical natural codes with singlefold validity, addition, negation and an induced height satisfying the hypotheses of Theorem~\ref{thm:abstract-height}. Consequently the same three conclusions hold on the quotient.
\end{corollary}
\begin{proof}
Insert the finitely many elements of $T$ as constants. For each $Q\in G$ compute the codes of all $Q+t$, $t\in T$, and select their lexicographically least code, denoted $\pi(Q)$. Lexicographic comparison is singlefold by Lemma~\ref{lem:compiler}. A fixed chain of such comparisons computes the least code with unique internal coordinates. Thus $\pi$ is a singlefold function, constant precisely on cosets of $T$. A code is canonical if it is valid and equals its own computed canonical representative. The graph of quotient addition on these codes is obtained by computing $Q+R$ in $G$ and then $\pi(Q+R)$; quotient negation is computed similarly. Their outputs and all witnesses are unique. The canonical representative of a coset of index $j$ is $jP+t_j$ for some $t_j\in T$. Evaluating $H$ on that code gives a singlefold height on the quotient, and the assumed uniform estimate holds for it. Since $G/T$ is infinite cyclic, Theorem~\ref{thm:abstract-height} applies.
\end{proof}
In particular every rank-one elliptic curve over $\Q$ admits this construction: use its finite torsion subgroup and the height described in the proof of Corollary~\ref{cor:positive-rank} below, with $p(j)=2\widehat{h}(P)j^2$. The canonicalization uses a finite list of codes rather than a coordinate ordering adapted to a particular curve. The torsion-free choice in Section~\ref{sec:height} remains convenient for the short direct proof.

\subsection{Canonical cyclic quotients in arbitrary positive rank}\label{sec:rank-quotient}
Suppose now that the group has rank greater than one. Passing to the quotient by a complementary subgroup leaves infinitely many points in each coset. Which point should stand for the coset, and how are we to recognize that point by means of finitely many polynomial constraints? It is natural to choose the point of least height, breaking ties by lexicographic order\index{lexicographic order} on the codes; the difficulty is to avoid an unbounded minimization. The key observation is that the bounded-error estimate for a logarithmic height with quadratic main term confines all least points to a region of uniformly bounded transverse width. This not only enables us to recognize them but also bounds the corrections required to perform group operations on them.

\begin{theorem}[Canonical cyclic quotient by height minimization]\label{thm:rank-quotient}\index{canonical cyclic quotient}
Let $G$ be a finitely generated Abelian group of rank $d\geq1$ with identity $O$ and torsion subgroup $T$. Assume unique natural codes of fixed length and singlefold polynomial systems for validity, addition, negation and a positive integer height $H:G\to\N_{>0}$. On $V=(G/T)\otimes_{\Z}\mathbb{R}$ suppose that there is a positive definite quadratic form $q$ such that
\begin{equation}\label{eq:rank-height}
 |\log H(R)-q(\bar R)|\leq C\qquad(R\in G),
\end{equation}
where $C\geq0$ is fixed and $\bar R$ denotes the image in $V$. Fix lifts $P_1,\ldots,P_d$ of an integral basis of $G/T$, and put $\mathcal{K}=\Z P_2+\cdots+\Z P_d+T$. Write $[P_1]$ for the class of $P_1$ modulo $\mathcal{K}$. Then $G/\mathcal{K}\cong\Z$ admits canonical natural codes with singlefold validity, addition, negation and a positive integer height $H_*$ satisfying
\[
 \log H_*(j[P_1])=aj^2+O(1),\qquad a>0.
\]
These systems are constructed from the supplied systems using finitely many fixed translations and elementary polynomial constraints. No membership test for $\mathcal{K}$ on arbitrary group elements is assumed.
\end{theorem}
\begin{proof}
Let $\operatorname{code}(R)$ denote the unique natural code of $R$. Write $\operatorname{key}(R)=(H(R),\operatorname{code}(R))$, ordered lexicographically.\index{key} In each coset choose the point of least key, and denote the representative of the coset of $jP_1$ by $R_j$. Such a representative exists: heights are positive integers, and there are only finitely many points below any given height bound by \eqref{eq:rank-height}, positive definiteness, discreteness of the free lattice and finiteness of $T$. The key distinguishes all points, even when heights tie.

\smallskip\noindent\emph{Uniform transverse bounds.} Every element is uniquely of the form $R=jP_1+\sum_{i=2}^d m_iP_i+\tau$, with $j,m_i\in\Z$ and $\tau\in T$. These are proof coordinates, not additional coded operations. Let $\langle\sqdot,\sqdot\rangle$ be the bilinear form of $q$, so $q(v)=\langle v,v\rangle$. If $H(R)\leq H(R\pm P_i)$ for $2\leq i\leq d$, expansion of \eqref{eq:rank-height} gives
\begin{equation}\label{eq:rank-local}
 |\langle\bar R,\bar P_i\rangle|\leq \tfrac12q(\bar P_i)+C.
\end{equation}
Indeed $q(\bar R\pm\bar P_i)-q(\bar R)\geq-2C$. In particular \eqref{eq:rank-local} holds for every $R_j$. For $d>1$ introduce the positive definite transverse Gram matrix and column vector
\[
 M=(\langle\bar P_i,\bar P_k\rangle)_{i,k=2}^d,\qquad
 h=(\langle\bar P_i,\bar P_1\rangle)_{i=2}^d,
 \qquad \beta=M^{-1}h.
\]
Write $m_i(R)$ for the integer coefficient of $P_i$ in the displayed expression for $R$. For a fixed coset the transverse coordinate vector at the minimum of the real quadratic form is $-j\beta$. Thus it is useful to measure transverse displacement\index{transverse coordinates} by $z(R)=(m_i)_{i=2}^d+j\beta$, rather than by the integers $m_i$ themselves. The vector of inner products in \eqref{eq:rank-local} equals $Mz(R)$; on multiplying by $M^{-1}$ and taking absolute values, we get the following bounds for any point that passes the local comparisons:
\begin{equation}\label{eq:rank-box-bound}
 |z_i(R)|\leq B_i:=\sum_{k=2}^d |(M^{-1})_{ik}|(M_{kk}/2+C).
\end{equation}
Choose fixed integers $N_i\geq\max(1,B_i)$. For each fixed positive integer $l$, form the finite list of distinct elements
\begin{equation}\label{eq:rank-boxes}
 F_l=\left\{\sum_{i=2}^d a_iP_i+\tau:\ a_i\in\Z,\ |a_i|\leq lN_i,\ \tau\in T\right\}.
\end{equation}
For $d=1$ all transverse components and matrices are omitted, the corresponding quadratic and bilinear terms below are zero, and $F_l=T$. The constants in the list can be inserted directly as natural codes. For $d>1$ its size is $|T|\prod_{i=2}^d(2lN_i+1)$.

\smallskip\noindent\emph{Recognizing the canonical domain.} A valid point $R$ is canonical exactly when
\begin{equation}\label{eq:rank-test}
 \operatorname{key}(R)\leq_{\rm lex}\operatorname{key}(R+L)\qquad(L\in F_2).
\end{equation}
A least point plainly passes this test. Conversely, $F_2$ contains the shifts $\pm P_i$, so the test implies \eqref{eq:rank-box-bound}. If $R$ belongs to coset $j$, then $z(R_j)-z(R)$ has components of absolute value at most $2N_i$. Thus the finite list $F_2$ contains the difference $R_j-R$ itself, irrespective of its torsion component. Indeed the displacement of two points in the same coset involves only the difference of their transverse coordinates and a torsion element; the term $j\beta$ cancels. One of the finite comparisons therefore compares $R$ with $R_j$ itself. Since the least key in the coset is that of $R_j$, the two keys must be equal, whence $R=R_j$. For $d=1$ the same argument reduces to comparing all torsion translates.

For the polynomial specification of this test, use the singlefold lexicographic comparison of Lemma~\ref{lem:compiler}. All list entries in \eqref{eq:rank-test} are fixed translations; every addition and height calculation is a singlefold call, and all comparison slacks and inactive branch coordinates are uniquely fixed by Lemma~\ref{lem:compiler}.

\smallskip\noindent\emph{Operations on the domain.} Addition in $G$ does not in general preserve the chosen representatives. We need not, however, search the entire coset of a sum. If $R=R_i$ and $S=R_j$, then $R_{i+j}-(R+S)$ belongs to $\mathcal{K}$. Its free transverse coordinate vector is $z(R_{i+j})-z(R)-z(S)$, each component of which has absolute value at most $3N_k$. Similarly $R_{-i}-(-R)$ has transverse components bounded by $2N_k$. Consequently
\begin{align}
 R\oplus S&=\underset{Q\in R+S+F_3}{\arg\min}\ \operatorname{key}(Q),\label{eq:rank-add}\\
 \ominus R&=\underset{Q\in -R+F_2}{\arg\min}\ \operatorname{key}(Q).\label{eq:rank-neg}
\end{align}
Here $R+S+F_3$ and $-R+F_2$ denote the corresponding translates of the finite sets; the minimization selects a point, not a translation. Every candidate lies in the required coset, and the globally least point is among them. Thus these are exactly the representatives of the sum and the opposite. Require canonical input codes by \eqref{eq:rank-test}, compute all candidates, and select the least by a fixed chain of deterministic comparisons. Each operation is uniquely determined together with all witnesses; no arbitrary existential choice of a translation is necessary. The zero and generator of this coded group are $R_0$ and $R_1$, respectively, not necessarily the original $O$ and $P_1$. They may be inserted as fixed codes (or found by finite minimization using the bounds $|m_k(R_0)|\leq N_k$ and $|m_k(R_1)|\leq N_k+|\beta_k|$).

\smallskip\noindent\emph{Height on the quotient.} Define $H_*(j[P_1])=H(R_j)$. Let $U$ be the real span of $\bar P_2,\ldots,\bar P_d$ and let $p_U$ be the orthogonal projection onto $U$ in the inner product of $q$. Put $v_1=\bar P_1-p_U(\bar P_1)$ and $a=q(v_1)>0$. Orthogonal decomposition gives\index{Schur complement}
\[
 q(\bar R_j)=aj^2+z(R_j)^t Mz(R_j),\qquad
 a=q(\bar P_1)-h^tM^{-1}h.
\]
The transverse contribution in the first expression is nonnegative and at most $C_M=\sum_{i,k=2}^d|M_{ik}|N_iN_k$ (zero in rank one). Thus $|\log H_*(j[P_1])-aj^2|\leq C+C_M$ for all signed indices. The graph of $H_*$ uses the canonical-domain test and the supplied height graph. This completes the construction.
\end{proof}

To turn the argument into a polynomial specification for a particular group, we need an integral basis of the free quotient and fixed upper bounds in \eqref{eq:rank-box-bound}; we are not claiming a uniform algorithm to extract these data from arbitrary input systems. Certified upper estimates suffice instead of exact real-number comparisons. Once these data and the finite torsion list are given, the polynomial construction is finite. Notice also that the theorem recognizes canonical points and performs operations on them; it does not canonically reduce every arbitrary point by a fixed finite list of translations. For example, $mP_2$ may require an unbounded translation to reach the chosen representative of the zero coset.

\begin{corollary}[Elliptic curves of positive rank]\label{cor:positive-rank}
For every fixed elliptic curve $E/\Q$ of positive Mordell--Weil rank, its naturally coded rational-point group and the full integer $x$-height yield a singlefold coded cyclic group satisfying the hypotheses of Theorem~\ref{thm:abstract-height}. In particular they yield a total function with a singlefold polynomial graph and logarithmic-square growth by the relative construction of that theorem.
\end{corollary}
\begin{proof}
Use a fixed nonsingular short Weierstrass model over $\Q$ and the rational codes and group law in Appendix~\ref{app:codes}, with its coefficients in place of those of \eqref{eq:curve}. Use height $1$ at $\OO$ and $H(x(Q))$ at an affine point $Q$. Mordell--Weil gives finite generation and the canonical height gives the positive definite quadratic form $2\widehat{h}$ on the free quotient, with a uniform difference from the logarithmic $x$-height. Theorem~\ref{thm:rank-quotient} supplies the hypotheses of Theorem~\ref{thm:abstract-height} with a quadratic main term. Its symmetrization, normalization and fixed amplification apply even if the quotient representatives do not have symmetric heights.
\end{proof}

For a rank-two illustration, write $q(j,m)=\alpha j^2+2\lambda jm+\gamma m^2$, with $\gamma>0$ and $\alpha\gamma>\lambda^2$. The entire transverse test is controlled by $|m+j\lambda/\gamma|\leq 1/2+C/\gamma$; the induced coefficient is $a=\alpha-\lambda^2/\gamma>0$. A concrete example is Cremona's curve 389a1, $y^2+y=x^3+x^2-2x$, of rank two with trivial torsion and basis $(0,0),(1,0)$ according to the published curve data \cite{LMFDB389}. It becomes $Y^2=X^3-3024X+46224$ under $X=36x+12$, $Y=216y+108$. Taking the subgroup generated by the second basis point as the complementary subgroup gives the stated cyclic construction even though the full rational-point group is not cyclic. This illustration uses published Mordell--Weil data; it does not supply certified numerical bounds for the height error and matrix entries, or an explicit list of translations for this curve.

\paragraph{Why the hypotheses matter.} Comparisons only with the coordinate generators are not sufficient to recognize a global minimum. On $G=\Z^3$, consider the form $q(j,m,n)=j^2+10m^2+18mn+10n^2$ and the illustrative height $H(j,m,n)=2^{q(j,m,n)}$. For this height, the point $(0,1,-1)$ has exponent $2$, below all four coordinate neighbours, yet the origin has exponent zero. The larger list $F_2$ is needed in general. Nor can the bounded error in \eqref{eq:rank-height} simply be replaced by an asymptotic equivalent. On $\Z^2$, for $f(j)=\lfloor\sqrt{|j|}\rfloor$ and $H(j,m)=2^{j^2+(m-f(j))^2}$ one has $\log H(j,m)\sim(\log2)(j^2+m^2)$ as $j^2+m^2\to\infty$, but the minimal representatives are $(j,f(j))$. Their sum can require the unbounded transverse correction $f(2j)-2f(j)$. These examples concern the geometry of the construction, not a claim of singlefold polynomial definability for the illustrative heights before exponentiation has been obtained.

\section{Growth envelopes and transfer of witness multiplicities}\label{sec:transfer}
The role played by the growth function in Proposition~\ref{prop:pow} need not be filled by a total function with a univocal specification. A wider class of polynomial systems can serve the same purpose. We shall construct power graphs by scale equations and Pell constraints applied to the supplied systems, accounting for all their solving tuples, including those yielding insufficient cutoffs. Simply adjoining an independent power graph would not perform this decoding.

\subsection{Transfer under a polynomial bound on the logarithm}
For a given input, distinct solving tuples of the supplied system may give the same cutoff, or they may give different cutoffs. The latter possibility brings into relief the two roles of growth already encountered in \eqref{eq:g-bounds}: every output must be small enough to preclude extraneous Pell indices, whereas it suffices that one output be large enough to accommodate the intended ones.
\begin{definition}[Logarithmically bounded cutoff family]\label{def:log-envelope}\index{cutoff family!logarithmically bounded}
Let $B(t,v,\mathbf{z})=0$ be a fixed finite system of integer polynomial equations over $\N$. Write $\mathcal{B}_t=\{(v,\mathbf{z}):B(t,v,\mathbf{z})=0\}$. Suppose for some fixed real constants $C,\eta>0$, integer $D\geq1$ and integer threshold $T_0\geq2$, every $t\geq T_0$ satisfies
\begin{align}
 (v,\mathbf{z})\in\mathcal{B}_t&\ \Longrightarrow\ v\geq1,\quad \log v\leq Ct^D,\label{eq:log-upper}\\
 &\exists(v,\mathbf{z})\in\mathcal{B}_t:\quad \log v>\eta(\log t)^2.\label{eq:log-lower}
\end{align}
All internal witnesses of $B$ are included in $\mathbf{z}$. The bounds are numerical assumptions, not additional constraints in the system. We call the family finite-fibre\index{finite-fibre family}\index{k-bounded family@$k$-bounded family} if $\mathcal{B}_t$ is finite for every $t\geq T_0$, and $k$-bounded if $\#\mathcal{B}_t\leq k$ for every $t\geq T_0$, where $k$ is a fixed positive integer.
\end{definition}

\paragraph{Common lower-cutoff estimate.}\index{common lower-cutoff estimate} Let $b,u,n,t\in\N$ with $b,u,t\geq2$, and let $L,\ell,\eta$ be positive real constants. Suppose
\[
 2bu+2\leq t^L,\qquad \log t\geq\ell n\log2,\qquad \ell>\frac{L}{\eta\log2}.
\]
By \eqref{eq:pell-bounds},
\begin{align*}
 \max\{\log y_{n+1}(bu+1),\log y_{n+1}(u)\}
 &\leq nL\log t\\
 &\leq\frac{L}{\ell\log2}(\log t)^2<\eta(\log t)^2.
\end{align*}
Thus any cutoff $v$ with $\log v>\eta(\log t)^2$ exceeds both required second coordinates. This estimate applies to each active Pell block in all four constructions below, independently of the choice of its scale.

\begin{theorem}[Cutoff transfer with exact multiplicities]\label{thm:transfer}\label{thm:log-upper}\index{cutoff transfer}\index{exact witness transfer}
For any family of Definition~\ref{def:log-envelope}, there is a finite polynomial system representing $c=b^n$ for $b\geq2$, $n,c\in\N$. For suitable fixed positive integers $\ell,t_0$, put
\begin{equation}\label{eq:log-transfer-scale}
 t=(c+2)^\ell(b+1)(n+2)+t_0,\qquad u=t^{D+1}+3(c+1)(n+1)+5.
\end{equation}
The complete witness fibre of the new system is empty for $c\ne b^n$. For $c=b^n$ projection onto the supplied cutoff tuple gives a bijection with the nonempty set
\begin{equation}\label{eq:log-transfer-count}
 \{(v,\mathbf{z})\in\mathcal{B}_t:v>\max(y_{n+1}(bu+1),y_{n+1}(u))\}.
\end{equation}
In particular finite-fibre or $k$-bounded families give finitefold or $k$-bounded power graphs. If $\log v>\eta(\log t)^2$ for every $(v,\mathbf{z})\in\mathcal{B}_t$, the bijection is with the entire fibre $\mathcal{B}_t$.
\end{theorem}
\begin{proof}
Choose $t_0\geq T_0$ with $t_0>C+2$. Introduce $t,u$ by the polynomial equalities \eqref{eq:log-transfer-scale}, impose $b\geq2$, call $B(t,v,\mathbf{z})=0$, and append the paired-Pell system \eqref{eq:general-pell}. The constants $D+1$ and $\ell$ are fixed exponents. For every possible cutoff tuple, not only the sufficiently large ones,
\[
 \log v\leq Ct^D<t^{D+1}<(u-1)\log(2u-1).
\]
Also $u\geq3(c+1)(n+1)$. Lemma~\ref{lem:exact-cutoff} therefore excludes every false candidate $c$ and, when $c=b^n$, gives exactly one completion for each member of \eqref{eq:log-transfer-count} and no completion for the others. All additional slacks, scale variables and Pell witnesses are uniquely determined given the cutoff tuple and the displayed arguments. Projection and unique completion are inverse operations. In particular, we have not identified two distinct witnesses merely because they yield the same cutoff.

To show that at least one cutoff is admissible, take any $\ell\geq1$ initially. For $t\geq2$, the scale gives $b<t$, $(c+1)(n+1)\leq t$ and $u\leq4t^{D+1}$. Thus $2bu+2\leq9t^{D+2}\leq t^{D+6}$. Set $L=D+6$, independently of $\ell$, and fix $\ell>L/(\eta\log 2)$. For a true instance $c=b^n$ we have $\log t\geq\ell n\log2$ so the common lower-cutoff estimate applies.
Every cutoff tuple with $\log v>\eta(\log t)^2$ is therefore admissible. Such a tuple exists by \eqref{eq:log-lower}. If all cutoff tuples satisfy the lower estimate, all of them survive the construction.
\end{proof}

The entire fibre may even be infinite while the representation of Theorem~\ref{thm:transfer} is finitefold. For example, append to the graph of $g$ from Proposition~\ref{prop:growth} a disjoint branch $v=1$ with a free local witness (restrict to sufficiently large $t$ so that $g(t)>1$). That branch makes the supplied fibre infinite but never passes an active Pell block; the transferred graph for $b\geq2$ is still singlefold. It would not be discarded by the common-cutoff construction on an entirely elementary conjunction, where there is no Pell block. Conversely the condition that every member satisfy the numerical lower estimate is sufficient, not necessary, for the full-fibre bijection on a particular instance: at $n=0$ any $v>1$ is admissible.

Could we dispense with the upper estimate altogether? Consider the following failure of this particular implementation. With $b=2$, $n=c=0$, $u=3$ and $v=36$, the seven variables in \eqref{eq:general-pell} can take the values $(X,Y,F,H,q,q',r)=(7,99,1,35,0,17,1)$, although $c\ne2^0$. The second Pell index is $3$, which the upper estimate in Lemma~\ref{lem:exact-cutoff} would exclude. This is a limitation of the prescribed-cutoff system, not a proof that every other transfer method must impose an upper bound.

\subsection{A single cutoff for a finite conjunction of power equations}\label{sec:joint-transfer}
The reduction stated in \cite[Theorem 5]{CCO2024}, after Matiyasevich, has the notable feature of relegating exponentiation in a whole conjunction of power equations to a sole literal. If we applied Theorem~\ref{thm:transfer} separately to each of the equations, the multiplicities coming from the several cutoff fibres could be multiplied together. Letting all Pell blocks share one cutoff avoids this multiplication of witnesses. The following result recovers the one-growth-constraint format under our hypotheses and counts distinct witnesses even when they yield identical cutoff values.

\begin{theorem}[Common-cutoff transfer]\label{thm:joint-transfer}\index{cutoff transfer!common-cutoff}
Let $B(t,v,\mathbf{z})=0$ satisfy Definition~\ref{def:log-envelope}, and fix a positive integer $m$. There is a finite integer polynomial system over $\N$, with one occurrence of the supplied system $B$, representing
\begin{equation}\label{eq:joint-powers}
 \bigwedge_{i=1}^m b_i^{n_i}=c_i,
 \qquad b_i,n_i,c_i\in\N.
\end{equation}
For suitable fixed integers $\ell,t_0\geq1$, set
\begin{equation}\label{eq:joint-scale}
 \begin{aligned}
 t&=t_0+\sum_{i=1}^m(c_i+2)^\ell(b_i+1)(n_i+2),\\
 u&=t^{D+1}+3\sum_{i=1}^m(c_i+1)(n_i+1)+5.
 \end{aligned}
\end{equation}
For the fixed displayed arguments let $I=\{i:b_i\geq2\}$ and define the numerical threshold
\begin{equation}\label{eq:joint-threshold}
 \Lambda=\max\bigl(\{0\}\cup\{y_{n_i+1}(b_i u+1),y_{n_i+1}(u):i\in I\}\bigr).
\end{equation}
If any power equality is false, the system has no complete witness. If all are true, projection onto $(v,\mathbf{z})$ is a bijection from its complete witness fibre onto the nonempty set
\begin{equation}\label{eq:joint-count}
 \{(v,\mathbf{z})\in\mathcal{B}_t:v>\Lambda\}.
\end{equation}
If $\log v>\eta(\log t)^2$ for every $(v,\mathbf{z})\in\mathcal{B}_t$, the bijection is with the entire fibre. In particular, a finite-fibre or $k$-bounded family yields a finitefold or $k$-bounded representation of the conjunction, with no factor $k^m$.
\end{theorem}
\begin{proof}
Choose $t_0\geq T_0$ with $t_0>C+2$. Put $L=D+6$ and choose a fixed integer $\ell>\max\{1,L/(\eta\log2)\}$. Introduce $t,u$ by \eqref{eq:joint-scale}, and call $B(t,v,\mathbf{z})=0$ once. For each $i$ independently impose one of the following disjoint cases:
\[
 \begin{array}{ll}
 b_i=0,\ n_i=0,\ c_i=1; & b_i=0,\ n_i>0,\ c_i=0;\\
 b_i=1,\ c_i=1; & b_i\geq2,\ \text{system \eqref{eq:general-pell} at }(b_i,n_i,c_i,u,v).
 \end{array}
\]
Apply the selector construction of Lemma~\ref{lem:compiler} with separate local variables, including slacks, for every case. The shared variables $t,u,v,\mathbf{z}$ are not zeroed. The number of blocks and their sizes are fixed. Neither $I$ nor $\Lambda$ has to be computed by a polynomial system: they occur only in the description of the fibre, not in the equations.

For every cutoff tuple the upper estimate gives
\[
 \log v\leq Ct^D<t^{D+1}<(u-1)\log(2u-1).
\]
Moreover $u\geq3(c_i+1)(n_i+1)$ for each $i$. For any block with $b_i\geq2$, Lemma~\ref{lem:exact-cutoff} therefore implies its power equality and fixes its complete Pell tuple precisely when both required coordinates are smaller than $v$. The other blocks simply impose the elementary equalities for bases zero and one. Thus every false conjunction has no witnesses. For a true conjunction, a supplied tuple has a unique completion if $v>\Lambda$ and no completion otherwise: all the selectors, inactive variables, active local witnesses, slacks and scales are fixed once the displayed arguments and the cutoff tuple are fixed. Conversely each complete solution determines a member of \eqref{eq:joint-count}. This proves the asserted bijection, including multiplicities of identical cutoff values with distinct internal witnesses.

It remains to ensure that some cutoff works simultaneously for all the equations. Observe that $b_i<t$ and $\sum_i(c_i+1)(n_i+1)\leq t$, so $u\leq4t^{D+1}$ and $2b_i u+2\leq9t^{D+2}\leq t^L$ for $i\in I$. These bounds do not depend on $m$ or on $\ell$. On a true conjunction, for every $i\in I$ we have $\log t\geq\ell n_i\log2$ since $c_i=b_i^{n_i}$. Hence the common lower-cutoff estimate applies to every active block.
The same lower estimate suffices for every block, so one and the same member of $\mathcal{B}_t$ meeting that estimate gives the required cutoff. If $I$ is empty, $\Lambda=0$ and every cutoff tuple is admissible, by positivity; the call is retained even in that case to keep the fibre formula uniform. Finally if every member satisfies the lower estimate, all members pass the common cutoff. This proves all assertions.
\end{proof}

The theorem also applies when the displayed arguments of different power equations reuse variables or are fixed natural constants: substitute the corresponding arguments in the polynomial system. The statement counts witnesses after all the displayed parameters are fixed; it does not assert uniqueness of an unconstrained existential choice of an argument. The difference between separate and shared calls already appears for a family obtained from a sufficiently large cutoff graph by adjoining a free Boolean witness $e$ with $e(e-1)=0$. If both choices have the same admissible cutoff, separate calls on $m$ non-elementary true power equations have $2^m$ witnesses, whereas the common call has exactly two.

\subsection{Retaining a scale-dependent cutoff}
In connection with Robinson's growth condition, it is convenient to have a designated scale among the witnesses of the supplied system. This suggests using the scale itself in the Pell parameter, rather than bounding its effect via $\log v$. Both constructions give a bijection, but in general the subsets of the supplied fibre singled out by them differ.

\begin{definition}[A polynomial cutoff family]\label{def:envelope}\index{cutoff family!polynomial}
A polynomial cutoff family is a fixed finite system of integer polynomial equations $B(t,s,v,\mathbf{z})=0$ (the notation permits a single equation), with all variables natural. The adjective polynomial refers to the supplied system, not to the growth of its cutoffs. For each $t$, write
\[
 \mathcal{B}_t=\{(s,v,\mathbf{z}):B(t,s,v,\mathbf{z})=0\}.
\]
Assume that for some fixed integer threshold $T_0\geq2$, real $A\geq1$, integer $d\geq1$ and real $\eta>0$ the following conditions hold for every $t\geq T_0$:
\begin{align}
 (s,v,\mathbf{z})\in\mathcal{B}_t&\ \Longrightarrow\ 2\leq s\leq A t^d,\quad 1\leq v\leq s^s,\label{eq:envelope-upper}\\
 &\exists (s,v,\mathbf{z})\in\mathcal{B}_t:\quad \log v>\eta(\log t)^2.\label{eq:envelope-lower}
\end{align}
The coordinates $\mathbf{z}$ include all witnesses of the supplied system; the variable $s$ is a scale and $v$ a cutoff. These estimates are hypotheses about the solutions, not predicates to be implemented with variable exponentiation. The conventions concerning finite-fibre and $k$-bounded families are those of Definition~\ref{def:log-envelope}, applied to the complete fibres $(s,v,\mathbf{z})$ for $t\geq T_0$.
\end{definition}
The upper estimate gives $\log v\leq s\log s\leq A^2t^{2d}$, so the preceding theorem applies with $s$ retained as an additional witness. The following construction instead lets the Pell scale depend on $s$.

\begin{corollary}[Scale-dependent cutoff transfer]\label{cor:scale-transfer}\index{cutoff transfer!scale-dependent}
Given a polynomial cutoff family, there is a finite polynomial system for $b^n=c$, $b\geq2$, with the following property. There are fixed positive integers $\ell,t_0$ such that for fixed $b,n,c$ in this domain put
\begin{equation}\label{eq:transfer-scale}
 t=(c+2)^\ell(b+1)(n+2)+t_0,\qquad u_s=s+3(c+1)(n+1)+5.
\end{equation}
If $c\ne b^n$, the system has no witnesses. If $c=b^n$, its complete witness set is in bijection with
\begin{equation}\label{eq:transfer-count}
 \{(s,v,\mathbf{z})\in\mathcal{B}_t:
 v>\max(y_{n+1}(bu_s+1),y_{n+1}(u_s))\},
\end{equation}
which is nonempty. Thus finite-fibre families give finitefold power graphs, and $k$-bounded families give power graphs with between one and $k$ complete witnesses per true instance. If $\log v>\eta(\log t)^2$ for every $(s,v,\mathbf{z})\in\mathcal{B}_t$, the bijection on true instances is with the entire $\mathcal{B}_t$.
\end{corollary}
\begin{proof}
Let $L$ be a fixed positive integer so large that
\[
 2(A+6)t^{d+1}+2\leq t^L\qquad(t\geq2).
\]
Such an integer exists, for example by increasing $L$ until $2^{L-d-1}\geq2(A+6)+2$. Choose a positive integer $\ell>L/(\eta\log2)$ and $t_0\geq T_0$. Introduce $t,u$ by the polynomial equations \eqref{eq:transfer-scale} (with $u=u_s$), impose $b\geq2$, call $B(t,s,v,\mathbf{z})=0$, and append \eqref{eq:general-pell}. The constant exponent $\ell$ is fixed. Every complete solution comes from one member of $\mathcal{B}_t$. For it,
\[
 u\geq3(c+1)(n+1),\qquad v\leq s^s<(2u-1)^{u-1},
\]
since $2\leq s<u-1$. Lemma~\ref{lem:exact-cutoff} proves soundness and gives a unique completion precisely for \eqref{eq:transfer-count} on a true instance. The unique-completion argument of Theorem~\ref{thm:transfer} then proves the bijection.

It remains to find an admissible tuple when $c=b^n$. In \eqref{eq:transfer-scale} we have $b<t$, $(c+1)(n+1)\leq t$, and $s\leq A t^d$. Hence $u\leq A t^d+3t+5\leq(A+6)t^d$ and $2bu+2\leq t^L$. Together with $\log t\geq\ell n\log2$ on a true instance, the common lower-cutoff estimate applies. The lower bound \eqref{eq:envelope-lower} ensures a nonempty admissible subfibre, and if every member satisfies it the entire fibre survives.
\end{proof}

Both transfer constructions extend to bases zero and one. Use the four disjoint cases in the proof of Theorem~\ref{thm:joint-transfer}, with the respective transfer system in the active branch. By Lemma~\ref{lem:compiler}, each true elementary branch contributes exactly one tuple, so the upper multiplicity bound is unchanged when $k\geq1$. The fibre formula of the transfer theorem is asserted for $b\geq2$, not for the elementary branches. The substitution principle following Theorem~\ref{thm:mat} transfers these bounds, using the base-two instance, to every r.e. relation. No general implication from mere Robinson growth is asserted.

The common-cutoff construction also has a scale-dependent version. We state it separately to distinguish its admissibility threshold from that of Theorem~\ref{thm:joint-transfer}.
\begin{corollary}[Scale-dependent common-cutoff transfer]\label{cor:joint-scale-transfer}
Let $B(t,s,v,\mathbf{z})=0$ satisfy Definition~\ref{def:envelope}, and fix $m\geq1$. There is a finite polynomial system representing \eqref{eq:joint-powers}, with one occurrence of $B$. For suitable fixed positive integers $\ell,t_0$, put
\begin{equation}\label{eq:joint-dependent-scale}
 \begin{aligned}
 t&=t_0+\sum_{i=1}^m(c_i+2)^\ell(b_i+1)(n_i+2),\\
 u_s&=s+3\sum_{i=1}^m(c_i+1)(n_i+1)+5.
 \end{aligned}
\end{equation}
For each supplied tuple let $I=\{i:b_i\geq2\}$ and
\begin{equation}\label{eq:joint-dependent-threshold}
 \Lambda_s=\max\bigl(\{0\}\cup\{y_{n_i+1}(b_i u_s+1),y_{n_i+1}(u_s):i\in I\}\bigr).
\end{equation}
The system has no witnesses for a false conjunction. For a true conjunction its complete witness fibre projects bijectively onto the nonempty set
\begin{equation}\label{eq:joint-dependent-count}
 \{(s,v,\mathbf{z})\in\mathcal{B}_t:v>\Lambda_s\}.
\end{equation}
A finite-fibre (respectively $k$-bounded) family therefore yields a finitefold (respectively $k$-bounded) representation. If $\log v>\eta(\log t)^2$ for every $(s,v,\mathbf{z})\in\mathcal{B}_t$, the projection is a bijection with the entire fibre.
\end{corollary}
\begin{proof}
Use one call to $B(t,s,v,\mathbf{z})=0$, fix the scales by \eqref{eq:joint-dependent-scale}, and impose the four disjoint cases of Theorem~\ref{thm:joint-transfer} for every block, with $u_s$ in place of $u$. For each supplied tuple, $2\leq s<u_s-1$, so
\[
 v\leq s^s<(2u_s-1)^{u_s-1},\qquad u_s\geq3(c_i+1)(n_i+1)\quad(1\leq i\leq m).
\]
These are the upper-cutoff and scale hypotheses of Lemma~\ref{lem:exact-cutoff} for every active block. The unique-completion argument of Theorem~\ref{thm:joint-transfer} therefore applies unchanged, excluding false conjunctions and giving precisely \eqref{eq:joint-dependent-count} on true ones.

To ensure the existence of a common cutoff, take $t_0\geq T_0$ from Definition~\ref{def:envelope}. Since $\sum_i(c_i+1)(n_i+1)\leq t$ and $b_i<t$, we have $u_s\leq(A+6)t^d$ and $2b_i u_s+2\leq t^L$, with $L$ chosen as in Corollary~\ref{cor:scale-transfer}, independently of $\ell$. Fix $\ell>L/(\eta\log2)$. For each $i\in I$ on a true conjunction, $\log t\geq\ell n_i\log2$, so the common lower-cutoff estimate applies to every active block. By \eqref{eq:envelope-lower}, there is a cutoff with $\log v>\eta(\log t)^2$; it passes every active block. If $I$ is empty, $\Lambda_s=0$ and every cutoff tuple is admissible. This proves nonemptiness and the full-fibre assertion, as well as the multiplicity bound.
\end{proof}

The four variants are collected below. Since their Pell scales differ, the admissible subfibres need not coincide. In the scale-independent versions a scale, if present, can be retained as part of the supplied tuple.
\begin{center}
\small
\begin{tabular}{@{}p{.30\linewidth}p{.34\linewidth}p{.29\linewidth}@{}}
\toprule
Construction & Upper hypothesis & Bound from a $k$-bounded family\\
\midrule
Theorem~\ref{thm:transfer} & $\log v\leq Ct^D$ & $k$ (one power)\\
Theorem~\ref{thm:joint-transfer} & $\log v\leq Ct^D$ & $k$ (all $m$ powers)\\
Corollary~\ref{cor:scale-transfer} & $s\leq At^d$, $v\leq s^s$ & $k$ (one power)\\
Corollary~\ref{cor:joint-scale-transfer} & $s\leq At^d$, $v\leq s^s$ & $k$ (all $m$ powers)\\
\bottomrule
\end{tabular}
\end{center}
The bounds in the last column count complete witnesses, not distinct outputs or computational resources.

\subsection{A criterion stated directly in terms of a dyadic relation}
\begin{corollary}[A logarithmic-square envelope]\label{cor:sparse-envelope}\index{logarithmic-square envelope}
Let $M(s,v)$ have a finitefold polynomial representation and satisfy
\[
 M(s,v)\ \Longrightarrow\ s\geq2,\quad 1\leq v\leq s^s.
\]
Suppose that for some fixed real $A\geq1$, integer $d\geq1$, real $\eta>0$ and every sufficiently large integer $t$ there exist $s,v$ with
\begin{equation}\label{eq:sparse-envelope}
 M(s,v),\qquad s\leq A t^d,\qquad \log v>\eta(\log t)^2.
\end{equation}
Then exponentiation and every recursively enumerable relation have finitefold polynomial representations.
\end{corollary}
\begin{proof}
Round $A$ up to an integer. Form $B(t,s,v,\mathbf{z})$ by adjoining $s\leq A t^d$ to the supplied representation of $M(s,v)$, counting the unique slack in $\mathbf{z}$. This is a cutoff family. For each fixed sufficiently large $t$ only finitely many $s$ are possible; for each such $s$ the hypothesis bounds $v$ by $s^s$, and each fixed pair has finitely many complete witnesses. Thus $\mathcal{B}_t$ is finite for each sufficiently large $t$, so the family is finite-fibre. Corollary~\ref{cor:joint-scale-transfer} with $m=1$ gives exponentiation for all natural bases. The substitution principle following Theorem~\ref{thm:mat} gives the assertion for r.e. relations.
\end{proof}

For instance, fix a sufficiently large integer $U\geq2$. The restricted graph $M(s,v)$ given by $s\geq U$ and $v=g(s)$ satisfies \eqref{eq:sparse-envelope}, with $s=t$ and $0<\eta<1/a$; imposing $s=t$ by an equality makes the family singlefold. A sparse-domain restriction\index{sparse domain} works as well: additionally require $s$ to be a perfect square. The square predicate has a unique natural root, so the resulting relation on pairs is singlefold. At each large $t$ take $s=t^2$ and $v=g(t^2)$, for which $\log v\sim4(\log t)^2/a$. Choosing $s=t^2$ by an equality in the family (rather than allowing every scale below a bound) makes the entire family singlefold and recovers the singlefold conclusion. This example distinguishes singlefold-ness of a binary relation from singlefold-ness of the scale-and-cutoff selection; the former alone would not force the latter.

\begin{corollary}[A relative one-growth-constraint format]\label{cor:one-growth-call}\index{one-growth-constraint format}
Let $M\subseteq\N^2$ satisfy the numerical hypotheses of Corollary~\ref{cor:sparse-envelope} (without initially assuming a polynomial representation of $M$). For each fixed positive integer $m$, write $\mathbf{b}=(b_1,\ldots,b_m)$, $\mathbf{n}=(n_1,\ldots,n_m)$ and $\mathbf{c}=(c_1,\ldots,c_m)$. There is an integer polynomial $D_m$ such that \eqref{eq:joint-powers} is equivalent to
\begin{equation}\label{eq:one-growth-call}
 \exists s,v,\vw\ \bigl(M(s,v)\ \wedge\ D_m(\mathbf{b},\mathbf{n},\mathbf{c},s,v,\vw)=0\bigr),
\end{equation}
with all quantified arguments natural. For fixed displayed parameters and a fixed pair $(s,v)$ belonging to $M$, there is at most one $\vw$ satisfying $D_m=0$; it exists exactly when the conjunction is true, $s\leq At^d$, and $v>\Lambda_s$, where $\Lambda_s$ is defined in \eqref{eq:joint-dependent-threshold}. Here $t$ is determined by the displayed parameters through \eqref{eq:joint-dependent-scale}, with suitable fixed positive integers $\ell,t_0$ chosen as in Corollary~\ref{cor:joint-scale-transfer}, and $A$ is rounded upward if needed. If $P_M(s,v,\vz)=0$ is a supplied finitefold polynomial representation of $M$, substitution for the single occurrence of $M$ yields a polynomial representation whose complete witness fibre, on a true conjunction, projects bijectively onto
\[
 \{(s,v,\vz):P_M(s,v,\vz)=0,\ s\leq At^d,\ v>\Lambda_s\}.
\]
Here $\vz$ includes every witness of the supplied representation; on a false conjunction the fibre is empty.
\end{corollary}
\begin{proof}
Round $A$ upward as before. For the scale-dependent construction use $M(s,v)$ and the single comparison $s\leq At^d$ in place of the call to $B$. The proof of Corollary~\ref{cor:joint-scale-transfer} uses only the numerical bounds on the supplied pairs until the system for $M$ is substituted; it therefore applies with $M$ left as an external predicate. The remaining equations, selectors and comparisons are polynomial and collapse to a single $D_m$ by a sum of squares. For fixed $(s,v)$ in $M$, all auxiliary variables in $D_m$ are fixed precisely when $s$ passes the bound and $v>\Lambda_s$ on a true conjunction. Upon substitution, the same unique-completion argument gives the displayed bijection, retaining every internal witness of $M$. Finiteness follows as in Corollary~\ref{cor:sparse-envelope}. Alternatively, Theorem~\ref{thm:joint-transfer} applies with $C=A^2$ and $D=2d$ and with $s$ retained as a witness. Its scale-independent threshold is \eqref{eq:joint-threshold}; the resulting admissible subfibre need not agree with \eqref{eq:joint-dependent-count}.
\end{proof}

Here the uniqueness in \eqref{eq:one-growth-call} is relative to the chosen pair, not a claim that there is only one such pair. If the selection of $(s,v)$ and its internal witnesses is singlefold, the full representation is singlefold. More generally, a singlefold exponential specification consisting of polynomial equations and a fixed conjunction of power equations can be transformed by this corollary: on a true input the original complete tuple, including the arguments of all the power equations, is fixed, and the only possible multiplicity is that of the admissible cutoff fibre. This is the one-growth-constraint form sought in \cite[Section 7]{CCO2024}, under the sufficient hypotheses stated here.

\subsection{Amplification of weaker logarithmic growth}\label{sec:amplification}
A slower-growing relation may still meet the preceding criteria after a fixed number of successive uses: iteration can amplify the powers of the logarithm. No quadratic-height sequence is needed for the following argument.
\begin{theorem}[Iterated logarithmic-power criterion]\label{thm:iteration}\index{iterated logarithmic-power criterion}
Suppose a fixed polynomial system $J(x,y,\mathbf{z})=0$ over $\N$ has finitely many complete $(y,\mathbf{z})$ for each sufficiently large input $x$. Suppose there exist real numbers $1<p\leq q<\infty$, $A_0,B_0>0$ and an integer $T\geq2$ such that, for all integers $x\geq T$,
\begin{enumerate}
\item some complete solution has $y\geq1$ and $\log y\geq A_0(\log x)^p$;
\item every complete solution with $y\geq1$ satisfies $\log y\leq B_0(\log x)^q$.
\end{enumerate}
Then exponentiation and all recursively enumerable relations have finitefold polynomial representations. If the total number of complete $(y,\mathbf{z})$ at each $x\geq T$ is at most a fixed positive integer $k$, these representations can be chosen with at most $k^r$ complete witnesses on every input, where $r$ is any fixed positive integer with $p^r\geq2$. In particular, if $J$ represents the graph of a function with unique complete witnesses at each input in the stated domain, these representations are singlefold.
\end{theorem}
\begin{proof}
Increase $T$ so that $A_0(\log x)^p\geq\log x$ for $x\geq T$, and so that the finiteness hypothesis holds for every input $x\geq T$. Fix $r$ with $p^r\geq2$, and introduce a chain of $r$ calls
\[
 x_0=t,\quad J(x_{i-1},x_i,\mathbf{z}_i)=0,\quad x_i\geq T
 \quad(1\leq i\leq r).
\]
The chain is subject to the lower threshold at every stage, lest an intermediate value fall outside the domain on which our estimates hold. For each $t\geq T$ there exists a chain choosing a lower-bound witness at each step, since such a choice is at least its input, hence at least $T$. Induction gives positive constants $a_i,b_i$ depending only on $i,A_0,B_0,p,q$, with
\[
 \log x_i\geq a_i(\log t)^{p^i}
 \quad\text{along such a chain},\qquad
 \log x_i\leq b_i(\log t)^{q^i}\quad\text{on every chain}.
\]
Indeed the lower and upper constants obey $a_i=A_0a_{i-1}^p$ and $b_i=B_0b_{i-1}^q$, with $a_0=b_0=1$. Choose $0<\eta<a_r$. If $p^r=2$ this choice gives the strict lower estimate even at the critical exponent; if $p^r>2$ increase the threshold until it holds. Thus some chain has $\log x_r>\eta(\log t)^2$ eventually, and every chain has $\log x_r\leq b_r(\log t)^{q^r}<t\log t$ eventually. Form the family of Definition~\ref{def:envelope} by $s=t$, $v=x_r$, keeping all intermediate values, inner witnesses and slacks. It is finite-fibre, because it has fixed depth and finite branching at every permitted input; a branching bound $k$ gives at most $k^r$ chains, counting all coordinates. Choose an integer $T_0\geq T$ large enough for both estimates to hold for every $t\geq T_0$. The family then satisfies both envelope inequalities, with upper scale bound $s=t$. Corollary~\ref{cor:joint-scale-transfer} with $m=1$ gives the power graph for all natural bases, preserving the bound $k^r$; the substitution principle following Theorem~\ref{thm:mat} transfers this bound to all r.e. relations. Singlefoldness here requires one complete $(y,\mathbf{z})$ at each input, not simply uniqueness of $\mathbf{z}$ for each pair $(x,y)$.
\end{proof}

For example, let $f$ be a function satisfying the singlefold-graph and upper-growth assumptions of Theorem~\ref{thm:iteration}. An eventual lower growth estimate $\log f(t)\geq C(\log t)^{1.1}$ with $C>0$ is sufficient: the criterion guarantees success with eight calls ($1.1^8>2$), whereas seven do not meet its exponent threshold ($1.1^7<2$). More generally, this theorem concerns a fixed composition length selected from the growth exponents, not a variable-length iteration whose verification would require a new encoding of computation. The lower bound $p>1$ is a sufficient hypothesis for this particular amplification argument, not a necessary condition for singlefold definability by other methods. For instance, if $\log f(t)\sim\log t\log\log t$, then for each fixed $r$ the logarithm of the $r$-fold iterate is asymptotic to $\log t(\log\log t)^r$, which is still $o((\log t)^2)$. This follows by induction from $\log(\log t(\log\log t)^i)\sim\log\log t$. The fixed-depth argument does not cover that growth rate. It also does not handle unrestricted sparse domains whose output need not be a permitted input to the next call.

\section{Further consequences}\label{sec:consequences}

\subsection{Undecidability even for quartic equations with at most one solution}
The preprint \cite[Corollary 1.2 and Section 8]{OpenAI2026} states an undecidability consequence with degree and number of variables both bounded. The degree bound can be chosen to be four by a classical transformation preserving uniqueness (cf. \cite{Jones1980,Jones1982} for universal polynomials).
\begin{corollary}\label{cor:quartic}\index{quartic equations!undecidability}
There exists an integer $M\geq1$ such that it is algorithmically undecidable whether a polynomial in $M$ variables with integer coefficients and total degree at most four has a natural solution, even under the promise that it has at most one such solution.
\end{corollary}
\begin{proof}
Choose a nondecidable recursively enumerable set $K\subseteq\N$ and a singlefold polynomial representation $P(a,\vw)=0$. Write $P=P^+-P^-$ with $P^+,P^-$ having nonnegative integer coefficients. Evaluate both polynomials by fixed finite arithmetic circuits using addition, multiplication, and constants. Introduce an auxiliary natural variable for every intermediate output and impose the corresponding equations, together with equality of the final outputs. Each equation is linear or quadratic, and for fixed $a,\vw$ every new variable is uniquely determined. The sum of squares of these equations gives a nonnegative integer polynomial $G(a,\vw,\vz)$ of total degree at most four with the same singlefold property. Specializing $a$ gives an effectively specified family in a fixed number of unknowns. Each equation has at most one natural solution and has one exactly for $a\in K$. A decision procedure under the promise would decide $K$.
\end{proof}

\begin{remark}[Varying only the right-hand side]
If $G$ is the nonnegative polynomial from the proof of Corollary~\ref{cor:quartic}, set $U(a,\vw,\vz)=(a+1)(1-G(a,\vw,\vz))$. For $c>0$, the equation $U=c$ forces $G=0$ and $a=c-1$. Thus $U$ is a fixed polynomial of degree at most five, each of its fibres $U^{-1}(c)$ in natural variables, for $c>0$, has at most one point, and it is undecidable whether such a fibre is nonempty. This uses the familiar positive-value device of Putnam\index{Putnam's positive-value device} (compare \cite[Remark 1]{CCO2024}).
\end{remark}

\section{Revisiting the 2024 programme on Diophantine singlefold specifications}\label{sec:cco-programme}
Let us reconsider the two objectives set out in \cite[Section 7, p.~602]{CCO2024}, following the earlier investigations \cite{CCFO2019,CCFO2021,CO2021}. For the one concerned with replacing the sole exponential literal by a suitable growth relation, Theorem~\ref{thm:joint-transfer} and Corollary~\ref{cor:one-growth-call} give a relative solution with exact multiplicity control. For the one concerned with exhibiting such a relation, we have the graph of Proposition~\ref{prop:growth}. Curiously, that graph does not meet the lower growth requirement advocated in 2024. To appreciate the difference, let us take a closer look at the two classical conditions.

\subsection{Two growth requirements and a different route}\label{sec:criteria}
The first growth requirement on a relation $M\subseteq\N^2$ is the one proposed by Julia Robinson\index{Robinson's growth condition}, tagged $\boldsymbol{(\dagger)}$ in \cite{CCO2024}:
\begin{equation}\label{eq:dagger}
 \begin{gathered}
 M(u,v)\ \Longrightarrow\ u>1\ \text{and}\ v\leq u^u,\\
 \forall h\in\N\ \exists u,v\in\N\ \bigl(M(u,v)\ \wedge\ u^h<v\bigr).
 \end{gathered}
\end{equation}
Davis's reduction \cite{Davis1963}, recalled in \cite[Theorem 4]{CCO2024}, permits a relation satisfying this condition to replace exponentiation in a Diophantine specification, but does not in general preserve the uniqueness of witnesses.

The second growth requirement\index{Matiyasevich's growth condition}, denoted $\boldsymbol{(\ddagger)}$ in \cite[Section 7]{CCO2024}, is attributed there to Matiyasevich \cite{Mat2010}. It asks for fixed integers $\alpha>1$, $\beta,\gamma\geq0$, and $\delta>0$ such that
\begin{equation}\label{eq:ddagger}
 \forall w\in\N_{>0}\ \exists u,v\in\N\quad
 M(u,v)\ \wedge\ u<\gamma w^\beta\ \wedge\ v>\delta\alpha^w.
\end{equation}
This condition motivated the search in \cite{CCO2024} for eight relations associated with quaternary quartic equations. In Matiyasevich's own formulation \cite[Sec.~4, condition (15$''$), p.~86 of the original publication]{Mat2010}, it is imposed upon an equation whose complete solution set is finite once $u$ is fixed. By iterating the growth relation a fixed number of times, he bounds the witnesses of an already available representation. This technique does not call for Robinson's upper bound. Our route is a different one: we use both a lower and an upper bound to construct a power specification directly by exact selection of Pell indices.

\begin{proposition}\label{prop:criteria}\index{growth conditions, separation of}
For any sufficiently large fixed integer $U\geq2$, the relation
\[
 M_g(u,v)\quad\Longleftrightarrow\quad u\geq U\ \wedge\ v=g(u)
\]
has a singlefold polynomial representation and satisfies \eqref{eq:dagger}, but does not satisfy \eqref{eq:ddagger}.
\end{proposition}
\begin{proof}
That the restricted graph has a singlefold system follows from Proposition~\ref{prop:growth} and the unique slack for the threshold inequality. For sufficiently large $u$, $\log g(u)<u\log u$, while $\log g(u)/\log u$ tends to infinity. These two statements give \eqref{eq:dagger} after fixing $U$. Suppose \eqref{eq:ddagger} holds for $M_g$. Necessarily $\gamma>0$. Since $\log g(u)=O((\log u)^2)$ for $u\geq U$, the bound $u<\gamma w^\beta$ would imply, uniformly over the permitted $u$,
\[
 \log v=\log g(u)=O\bigl((1+\log w)^2\bigr)=o(w).
\]
Here the constants can depend on $\beta,\gamma$, but not on $w$ or the permitted $u$. On the other hand, $v>\delta\alpha^w$ gives $\log v>w\log\alpha+\log\delta$. Since $\log\alpha>0$, this is a contradiction for sufficiently large $w$.
\end{proof}

The separation in Proposition~\ref{prop:criteria} is already implicit in \cite[Proposition 4.2, pp.~9--10]{OpenAI2026}: the original growth function also has a logarithm bounded above and below, for sufficiently large inputs, by positive constant multiples of $(\log u)^2$. The same restriction of its graph and the same asymptotic contradiction apply. Thus the separation is not a consequence peculiar to our choice of curve.

The restricted relation $M_g$ itself can be used in Proposition~\ref{prop:pow}: choose $u_0\geq U$ in addition to its other lower bounds and replace $v=g(u)$ by $M_g(u,v)$. The scale equation ensures $u\geq U$, so the proof and uniqueness argument are unchanged. Thus \eqref{eq:ddagger} is not necessary for a particular relation to serve as a cutoff in this construction; this is neither a converse to Matiyasevich's sufficient condition nor a claim that Robinson's condition alone suffices. Under the Robinson upper bound and a positivity guard on outputs, \eqref{eq:ddagger} implies \eqref{eq:sparse-envelope}, since $t\log\alpha+\log\delta$ eventually exceeds $\eta(\log t)^2$. The new criterion trades the upper-bound freedom of Matiyasevich's method for a weaker lower-growth demand and gives exact multiplicity information. Without the additional upper bound the hypotheses are not claimed to be ordered. The sparse square-domain restriction in Section~\ref{sec:transfer} also fails \eqref{eq:ddagger} by the same upper estimate as $M_g$.

\subsection{Relation to earlier work and scope of the contributions}\label{sec:prior}
\paragraph{Comparison with the September 2026 proof.} In \cite[Proposition 3.1, Lemma 4.1 and Proposition 4.2]{OpenAI2026}, the curve $y^2=x^3-25x$ has rank one and four torsion points. A fixed even scalar kills the torsion and enlarges the quadratic coefficient of the logarithmic heights; the generator and scalar are not numerically specified. Sections~\ref{sec:height} and \ref{sec:growth} replace those steps by a specified generator, a one-sided height estimate and an orientation test, while retaining the original index-square identity. The two Pell quotients (\cite[Lemma 5.1]{OpenAI2026} and \eqref{eq:quotient}) are alternative means of recovering a power; no quantitative superiority of the shifted quotient is asserted. The subsequent binary-containment and counter-machine constructions of \cite[Sections 6 and 7]{OpenAI2026} are replaced here by Theorem~\ref{thm:mat}. The original preprint already recognizes the classical reduction to binary exponentiation \cite[Section 1, p.~2]{OpenAI2026}; our simplification is to use it in the proof rather than to develop the computation encoding. This is a simplification at the level of arithmetic operations and proof modules, not a comparison of the degrees, numbers of variables or coefficient sizes of fully expanded representing polynomials.

The two principal extensions beyond that proof are the finite construction of canonical cyclic quotients in arbitrary positive rank (Theorem~\ref{thm:rank-quotient}) and the exact transfer of complete witness fibres (Theorems~\ref{thm:transfer} and \ref{thm:joint-transfer}, together with the scale-dependent variants). The latter uses Lemma~\ref{lem:exact-cutoff} to extend the mechanism of \cite[Proposition 5.2]{OpenAI2026} from a singlefold function to supplied systems with potentially many witnesses. These are relative results about the supplied representations and do not follow merely by invoking the general singlefold existence theorem.

Singlefold and finitefold definitions have also been studied in other witness domains. For example, Miller and Shlapentokh \cite{MillerShlapentokh2022} obtain singlefold definitions of r.e. (there called computably enumerable) sets of rational integers in certain characteristic-zero polynomial rings (see Section~4 of the arXiv preprint). The witnesses in those results range over the rings in question, not over $\N$. Pasten \cite{Pasten2022} studies analogues of the DPRM theorem for listably presented structures and positive existential interpretations. The Davis--Putnam--Robinson--Matiyasevich (DPRM) property, namely the Diophantine definability of recursively enumerable relations in the relevant setting, does not assert uniqueness of the complete tuple of polynomial witnesses. Our quotient construction instead assumes particular singlefold operations on natural codes and preserves that property by finite comparisons. These are related but distinct definability questions.

The use of elliptic curves of rank one in Diophantine definability long predates the September manuscript. Poonen \cite{Poonen2002} uses such curves and their denominators to study definability between rings of integers; Poonen and Shlapentokh \cite{PoonenShlapentokh2005} construct Diophantine models of the integers in certain large subrings of number fields. In particular, Shlapentokh \cite{Shlapentokh2012} explicitly studies multiplication of indices of elliptic multiples over large subrings. Those results concern existential definitions over the relevant rings themselves. Proposition~\ref{prop:indices}, in contrast, uses integer numerator--denominator codes and the full integer $x$-height, and counts every auxiliary natural coordinate. The specific contribution of Proposition~\ref{prop:indices} is the direct height-window formula, a polarization of the index-square identity of \cite[Proposition 4.2]{OpenAI2026}, rather than Diophantine treatment of elliptic indices as such.

The existence of singlefold graphs in Corollary~\ref{cor:hierarchy} is also an abstract consequence of Theorem~\ref{thm:main}, since the functions are computable: rational-point arithmetic and an increasing, unbounded height sequence allow the bracket to be found by search. Its additional content is the direct construction by a fixed number of height-window calls, independent of the power graph and of the general representation theorem. Likewise Corollary~\ref{cor:quartic} is a consequence of the singlefold theorem and the classical degree-four reduction, rather than a new degree-reduction technique. It gives no numerical bound for the number of unknowns and does not assert that the older explicit universal polynomials of \cite{Jones1980,Jones1982} are singlefold. The right-hand-side construction in the remark following that corollary uses the classical positive-value device already described in \cite[Remark 1]{CCO2024}.

The survey \cite[Sections 3--5]{CCFO2021} juxtaposes Robinson's growth criteria, several classical Pell reductions and the candidate quartics. Cantone and Omodeo \cite[Appendix A, p.~27]{CO2021} explicitly identify the need to refine Robinson's reduction so as to preserve finitefold-ness, and point to the classical Pell quotient as a possible route. These precedents explain the formulation of the two objectives in \cite{CCO2024}. The classical quotient is rederived in Appendix~\ref{app:quotient}.

The geometric part of the quotient argument also has a classical antecedent: finding a closest lattice point with respect to a positive definite quadratic form is tied to finite descriptions of Voronoi cells (cf. \cite{MicciancioVoulgaris2013}). In Theorem~\ref{thm:rank-quotient}, however, the logarithm of the height to be minimized differs from such a form by a bounded amount. The point of the transverse bounds is to retain a finite test for this perturbed height and to implement it, and the group operations on the chosen representatives, by singlefold constraints. No new algorithmic bound for closest-vector problems is claimed.

The classical single-exponential format and the question of replacing its exponential constraint by a growth relation are recorded in \cite[Theorem 5 and Section 7]{CCO2024}. The shared-call construction here retains the multiplicity bound of a single supplied fibre rather than multiplying the bounds for separate calls; the growth criteria also accommodate nonfunctional families and sparse domains. Iteration to bound witnesses has an explicit antecedent in \cite{Mat2010}; Theorem~\ref{thm:iteration} treats a slower logarithmic-power regime with direct Pell decoding and a bound on complete witness counts. These comparisons identify the extensions of the cited methods; they do not establish absolute bibliographic priority.

\subsection{The quartic-equation programme: status and corrections}
The eight candidate `rule-them-all' equations\index{quartic equations!rule-them-all} of \cite{CCFO2021,CCO2024} were conceived as sufficient routes to finitefold representability: if one of them has only finitely many integral solutions, a suitable relation of exponential growth becomes finitefold Diophantine, and finitefold representability of all r.e. sets follows. The implication goes in one direction only, so Theorem~\ref{thm:main} tells us nothing about the finiteness of their solution sets. Their study remains an arithmetical problem in its own right.

\paragraph{A correction to our earlier papers.} The relations associated with the quartics, denoted $\mathcal{J}_d$ in \cite{CCFO2021} and $\mathcal{M}_d$ in \cite[Section 8]{CCO2024}, combine membership in a sparse Pell sequence with a condition of the form ``$p\mid q$ and $p\geq2^{\nu_2(q)}$'', whereas their Diophantine specifications use $q=(2x+1)p$ with $x\in\N$. For positive $p,q$, the latter means that $p\mid q$ and that every power of $2$ dividing $q$ divides $p$, as correctly observed in the proof of \cite[Corollary 2]{CCO2024}; that is, $p\mid q$ and $\nu_2(p)=\nu_2(q)$. This is strictly stronger than $p\mid q$ together with $p\geq2^{\nu_2(q)}$, with which it is identified there. The simplest instance occurs for $d=2$. Let $\langle y_i\rangle_{i\in\N}=\langle0,2,12,70,408,\ldots\rangle$ be the increasing sequence of nonnegative integers $y$ for which $2y^2+1$ is a perfect square, and recall from \cite[Section 5]{CCFO2021} that $\mathcal{J}_2(p,q)$ means $q=y_{2^\ell}$, $p\geq2^{\ell+1}$ and $p\mid q$ for some $\ell$. Then $\mathcal{J}_2(17,408)$ holds, with $\ell=2$, since $17\geq8$ and $408=y_4=17\cdot24$; but $408/17$ is even, so the bi-implication for $\mathcal{J}_2$ displayed in \cite[p.~149]{CCFO2021} fails. The relation $\mathcal{M}_{67}$ of \cite[Corollary 2]{CCO2024} shows the same discrepancy at every level. Define the Pell sequence $(X_i,Y_i)$ by $X_i+Y_i\sqrt{67}=(48842+5967\sqrt{67})^i$ for $i\in\N$, and put $N=2^{2\ell+1}$ and $v=Y_N/351$ with $\ell>12$. The divisor $u=X_{N/2}$ of $v$ is odd, since $X_{2i}=2X_i^2-1$, and exceeds $2N$; yet $v/u=2Y_{N/2}/351$ is even, because $Y_N=2X_{N/2}Y_{N/2}$ and $351=Y_1/17$ divides every $Y_i$.

The remedy is to take the narrower relation, defined by $q=(2x+1)p$, as the object of study: it is what the specifications actually represent. It contains the pairs that witness growth and is contained in the original relation; for $\mathcal{J}_2$, for example, the pairs $(2^{\ell+1},y_{2^\ell})$ remain, since $\nu_2(y_{2^\ell})=\ell+1$. Hence every lower growth property established through such pairs, and every upper bound valid for all pairs of the original relation, passes to it. In addition, the constants stated in Theorems 7 and 8 of \cite{CCO2024} require adjustment at small $w$. For $w=1$, the restriction $\ell>12$ in the $d=67$ family forces $u\geq2^{2\ell+2}\geq2^{28}$, whereas the proposed upper bound is $u<2^7$; for $d=2$, $\ell>4$ forces $u\geq2^{\ell+1}\geq2^6$, whereas the proposed upper bound is $u<2^6$. Moreover, the proof of Theorem 7 bounds $y_{2^{2\ell+1}}$ from below, whereas the bound is required for $v=y_{2^{2\ell+1}}/351$. All this is repaired by enlarging the constants, as the following argument shows for $d=67$; the pairs used in those proofs, with $u=2^{\nu_2(v)}$, belong to the narrower relation.

Fix a sufficiently large integer $\ell_0\geq13$. Conditionally on a finitefold polynomial representation of the unary set $\{Y_{2^{2\ell+1}}/351:\ell\geq \ell_0\}$, retain the pairs $(s,v)$ with $v$ in that set and $v=(2x+1)s$. This relation is finitefold on fixed pairs, since $x$ is determined. The duplication identity gives $\nu_2(Y_{2^m})=m+1$ for $m\geq1$: the first doubling gives valuation $2$ (since $\nu_2(X_1)=1$ and $Y_1$ is odd), and each subsequent $X_{2^j}$ is odd. Thus, for $N=2^{2\ell+1}$, $s=2N$ is an allowed divisor of $v$, and every allowed divisor is at least $2N$. Writing $\rho=48842+5967\sqrt{67}$, the Pell formula gives $\log v=N\log\rho+O(1)$; for our sufficiently large fixed $\ell_0$ this is at most $s\log s$ for every allowed divisor. For any $w\geq1$ choose the smallest $N=2^{2\ell+1}$ with $\ell\geq \ell_0$ exceeding $C_1w$, with $C_1$ a sufficiently large fixed constant. Then $C_1w<N\leq C_2w$ for a fixed $C_2$, the progression having ratio four, so $s=2N$ is linearly bounded in $w$ while $v$ is exponentially large in $w$. Hence both the criterion of Corollary~\ref{cor:sparse-envelope} and $\boldsymbol{(\ddagger)}$ hold for this narrower relation. The argument corrects the growth comparison, conditionally on the indicated unary representation; it is not a proof of the implication from a quartic to that representation, nor of the finiteness of a quartic. We have not re-audited all eight quartic constructions.

\subsection{Questions left open by the method}
In view of Theorem~\ref{thm:main}, the existence of finitefold representations for all r.e. sets is settled. The methodological question is whether the explicit decoding and witness-transfer method extends to an arbitrary finitefold representation of a relation satisfying $\boldsymbol{(\dagger)}$ alone, using controlled calls to the supplied relation and no independent power graph or elliptic cutoff. This is a question about the stated decoding method, not an unrestricted existential question; the examples separating growth conditions do not settle it. For instance $v=s^{\lfloor\sqrt{\log s}\rfloor}$, restricted to $s\geq3$, satisfies $\boldsymbol{(\dagger)}$ but cannot meet a logarithmic-square lower bound when $s\leq At^d$, since then $\log v=O((\log t)^{3/2})$; this separates the hypotheses numerically, without showing that a representation is impossible. It would also be useful to find an independent source of polynomial logarithmic heights of degree greater than two, or a relative amplification that handles near-polynomial growth beyond any fixed number of calls. Finally, the finite-translation construction of Theorem~\ref{thm:rank-quotient} is existential for each fixed group with the stated height data; quantitative bounds and explicit polynomials for concrete higher-rank curves remain to be investigated.

\section*{Acknowledgements and AI disclosure}
\addcontentsline{toc}{section}{Acknowledgements and AI disclosure}
The authors take full responsibility for the content of this article, including its mathematical arguments and bibliographic references. Generative AI tools (OpenAI's ChatGPT (GPT-6.2) and Anthropic's Claude (Opus 5.5)) were used to assist with mathematical exploration, the development and checking of proofs, computational experiments, bibliographic research, and the drafting and revision of the manuscript. The elliptic-height method of \cite{OpenAI2026} is a separate source and is cited as such.

\appendix
\section{Polynomial implementation of rational-point arithmetic}\label{app:codes}
Up to this point we have used rational-point arithmetic as if it were available among the primitive constructs of our specifications. Let us justify that usage. Rather than expand a very large polynomial, we describe the individual constraints and show that they determine every variable introduced in their translation. Our treatment follows \cite[Section 3]{OpenAI2026}; see also the discussion of coprimality and disjunction in \cite[Section 2]{CCO2024}.

\subsection{Canonical rational codes}
We use the shifted modular-inverse code of \cite[Section 3.2, Lemma 3.2]{OpenAI2026}, with the explicit nonzero-divisor guard below.
\begin{lemma}\label{lem:rational}\index{rational code}
Every rational number has exactly one natural code $(p,n,d,t,s,q)$ satisfying
\begin{equation}\label{eq:code}
 pn=0,\qquad d=1+t+s,\qquad (p+n)t+d=1+qd.
\end{equation}
Its value is $(p-n)/d$. Rational additions, subtractions, multiplications, divisions by nonzero numbers, and comparisons can be implemented by polynomial systems on such codes, with uniquely determined outputs and internal witnesses wherever defined.
\end{lemma}
\begin{proof}
Set $A=p+n$. The first equation splits the signed numerator uniquely, the second asserts $d>0$ and $0\leq t<d$, and the third implies $At\equiv1\pmod d$. In particular $\gcd(A,d)=1$. The encoded fraction is therefore in lowest terms. Conversely, take any reduced fraction with positive denominator. Its numerator uniquely gives $p,n$, and the denominator gives $d$. For $d>1$ there is a unique inverse $t$ of $A$ modulo $d$ with $0\leq t<d$. The other two coordinates are forced by $s=d-1-t$ and $q=(At+d-1)/d$, which are nonnegative integers. If $d=1$, the equations force $t=s=q=0$. This includes the code of zero $(0,0,1,0,0,0)$; no other denominator is possible for zero.

Write $z_i=A_i/d_i$ for the rational value of the $i$th code, where $A_i=p_i-n_i$ is its signed numerator and $d_i$ its positive denominator. For instance, the rational equation $z_1+z_2=z_3$ is equivalent to
\[
 A_1d_2d_3+A_2d_1d_3=A_3d_1d_2.
\]
The code equations are imposed for all outputs. Multiplication and subtraction are handled by the same clearing of positive denominators. For division $z_1/z_2=z_3$, impose
\[
 A_1d_2d_3=A_2A_3d_1,\qquad p_2+n_2=1+\sigma,\quad \sigma\in\N.
\]
The second equation explicitly excludes a zero divisor and has a unique slack. Without it, the cleared equation could admit arbitrary outputs when both operands vanish. The resulting value, when the operation is defined, has exactly one code by the first part of the proof. Rational comparisons are integer polynomial comparisons after clearing positive denominators and so have unique slacks by Lemma~\ref{lem:compiler}. No freely chosen common denominator or unbounded B\'ezout coefficient is introduced.
\end{proof}

\subsection{Codes of points and the group law}
The case construction adapts \cite[Section 3.2, Lemma 3.3]{OpenAI2026} to the curve \eqref{eq:curve}.
A point code\index{point code} comprises two rational codes and a flag $o\in\N$ satisfying $o(o-1)=0$. When $o=1$, both rational coordinates are forced to zero, so the code represents $\OO$ uniquely. When $o=0$, impose the affine curve equation $y^2=x^3-74$ on those codes, after clearing denominators. The two requirements are gated by $o$ and $1-o$, respectively. The flag distinguishes the dummy affine coordinates of the point at infinity from ordinary affine points. Thus every point of $E(\Q)$ has exactly one code and every valid code represents a point.

To specify the sum of two points $Q_1,Q_2$ univocally, write $Q_i=(x_i,y_i)$ whenever $Q_i$ is affine ($i=1,2$). We distinguish the following mutually exclusive alternatives:
\begin{center}
\begin{tabular}{@{}p{.50\linewidth}p{.38\linewidth}@{}}
\toprule
Guard&Rule\\
\midrule
$Q_1=\OO$&Return $Q_2$\\
$Q_1\ne \OO$, $Q_2=\OO$&Return $Q_1$\\
Both affine, $x_1<x_2$&Secant formula\\
Both affine, $x_1>x_2$&Secant formula\\
Both affine, $x_1=x_2$, $y_1+y_2=0$&Return $\OO$\\
Both affine, $x_1=x_2$, $y_1=y_2>0$&Tangent formula\\
Both affine, $x_1=x_2$, $y_1=y_2<0$&Tangent formula\\
\bottomrule
\end{tabular}
\end{center}
At equal $x$-coordinates, the affine curve equation implies $y_1^2=y_2^2$, so the cases are exhaustive as well as disjoint. In the two secant cases set the rational slope $\lambda$ by
\[
 (x_2-x_1)\lambda=y_2-y_1.
\]
In the two tangent cases set it by
\[
 2y_1\lambda=3x_1^2.
\]
The slope coefficient is nonzero in either situation. (On the curve \eqref{eq:curve} no rational point has $y=0$, since $74$ is not a cube; the guards do not rely on this.) The output is the affine point $Q_3=(x_3,y_3)$ with coordinates
\[
 x_3=\lambda^2-x_1-x_2,\qquad y_3=\lambda(x_1-x_3)-y_1.
\]
These are the chord--tangent formulas on the chosen short Weierstrass model. For a fixed nonsingular model $y^2=x^3+Ax+B$ over $\Q$, use that equation in the validity test and replace the tangent equation by $2y_1\lambda=3x_1^2+A$. The guards, secant and output formulas are unchanged; fixed rational coefficients are cleared by multiplication by their fixed positive denominators. The conditions above, including guards, slope-code equations and output identities, become polynomial after denominators are cleared. The output must be a valid shared point code. Every case-local coordinate is gated, including its rational-code validity conditions, and is set to zero in an inactive case as in Lemma~\ref{lem:compiler}. In particular an inactive dummy code is not subject to a separate, ungated positive-denominator condition. The active slope, when required, has precisely one rational value and one code, and the same is true of the output. This proves singlefold-ness of point addition. Negation is the map $(x,y)\mapsto(x,-y)$ on affine points and fixes $\OO$, with its own disjoint polynomial cases if needed.

A fixed scalar multiplication $Q\mapsto mQ$ is computed by a fixed finite chain of additions. Its intermediate point codes and every slope, selector, and guard slack are unique. For $W(Q)$ in \eqref{eq:Wdef}, take the maximum of $p+n$ and $d$ in the rational code of the $x$-coordinate of $Q$. At $\OO$, the dummy code gives the same value 1 automatically. Ties are assigned to a single case in the maximum operation. Constants such as $\GG$ and $W(\GG)$ have fixed numerical codes; they can be inserted as integer coefficients, rather than variables with unspecified values.

\subsection{Accounting for the complete tuple}
\paragraph{Height graph.} In Proposition~\ref{prop:growth}, the outer unknowns are the two valid point codes of $Q,R$. One may compute the points $Q+\GG,R+\GG,2Q$ and the associated heights, and the heights of $Q,R$, with a fixed finite number of calls. The comparisons in \eqref{eq:bracket}, \eqref{eq:orient} (for $R$) and \eqref{eq:window} add only unique slacks. The two products $L,D$ and the final value $g(u)$ are fixed by equality. No unknown $k$ or $r$ and no unknown real height occurs anywhere. A call to the graph of $g$ therefore comprises these outer coordinates and all the internal tuples of the finitely many point and height calls, expanded as described above. Repeated calls can either be retained separately or share already determined outputs; both conventions preserve uniqueness.

\paragraph{Power graph.} In Proposition~\ref{prop:pow}, the outer unknowns are $u,v,X,Y,F,H,q,q',r$ and the inner tuple of the call $v=g(u)$, together with one slack for each of the three comparisons in \eqref{eq:cutoff} and \eqref{eq:division}. Any fresh names for polynomial expressions have their defining equalities. Substitution into the classical representation uses fresh variables. 

\paragraph{Transfer and other constructions.} Retain every coordinate of the supplied cutoff tuple (including the scale when present). Lemma~\ref{lem:exact-cutoff} fixes the seven Pell variables and three slacks of each active block; the selector convention of Lemma~\ref{lem:compiler} fixes every inactive local coordinate at zero. Scale equations and comparisons add only uniquely determined outputs and slacks. The shared cutoff tuple is retained only once in the common-cutoff constructions. Thus projection onto the supplied tuple neither identifies distinct witnesses nor hides extra multiplicities. The canonical-quotient and iteration constructions likewise retain the internal tuples of all constituent calls; see their proofs.

\section{An elementary proof of the quotient formula}\label{app:quotient}
For the sake of completeness, we derive the quotient formula of Lemma~\ref{lem:quotient}, retaining the classical constant 3 and explaining why the lower bound on $u$ can legitimately involve the proposed output $c$. Under its hypotheses, suppose $b,n,c,u\in\N$, with $b\geq1$ and $u\geq3(c+1)(n+1)$. For $i\geq1$ put $R_i(A)=y_{i+1}(A)/y_i(A)$, as in the proof of Lemma~\ref{lem:pell}. Recall that
\[
 R_1(A)=2A,\qquad R_{i+1}(A)=2A-\frac1{R_i(A)},\qquad
 2A-1\leq R_i(A)\leq2A.
\]
We claim that $R_i(bu+1)\geq bR_i(u)$. It holds for $i=1$, and the recurrence gives
\[
 R_{i+1}(bu+1)-bR_{i+1}(u)
 =2+\frac{b}{R_i(u)}-\frac1{R_i(bu+1)}>0,
\]
since $R_i(bu+1)\geq2(bu+1)-1>1$. The upper bound follows directly from the individual bounds:
\[
 \frac{R_i(bu+1)}{R_i(u)}\leq\frac{2bu+2}{2u-1}
 =b\left(1+\frac{b+2}{b(2u-1)}\right)
 \leq b\left(1+\frac3{2u-1}\right).
\]
Multiplying these inequalities for $i=1,\ldots,n$ and using $y_1(A)=1$, we obtain, with $T=y_{n+1}(bu+1)/y_{n+1}(u)$ and $z=3/(2u-1)$,
\[
 b^n\leq T\leq b^n(1+z)^n.
\]
The hypothesis on $u$ gives $2u-1>3n(c+1)$, so $nz<1$. The binomial theorem followed by the geometric series gives
\[
 (1+z)^n-1\leq\sum_{j=1}^n(nz)^j\leq\frac{nz}{1-nz}.
\]
Consequently,
\[
 0\leq T-b^n\leq\frac{3nb^n}{2u-1-3n}.
\]
We can treat the two implications together. If $c=\lfloor T\rfloor$, the lower bound implies $c\geq b^n$ since $b^n$ is integral. In the other direction we assume $c=b^n$, so the same bound holds. Under either assumption,
\[
 0\leq T-b^n\leq\frac{3nc}{2u-1-3n}<1,
\]
where the last inequality is precisely a consequence of $2u-1>3n(c+1)$. It follows that $\lfloor T\rfloor=b^n$. This proves both directions of \eqref{eq:quotient}. The case $n=0$ uses empty products and sums and gives $T=1$, so there is no exceptional case. All estimates in this appendix are used only in the proof; they introduce no auxiliary polynomial variables.

\section{A rank computation for the Mordell curve}\label{app:descent}
We verify the rank and torsion assertions for $E:y^2=x^3-74$ used in Section~\ref{sec:height}. The primes $7$ and $13$ are of good reduction. For $p=7$ the numbers of affine points above $x=0,1,\ldots,6$ are $0,2,2,2,2,2,2$; for $p=13$ the respective counts are $2,0,2,0,2,2,2,2,2,0,2,2,2$. Adding the point at infinity gives $13$ and $21$. The reduction map is injective on rational torsion at odd primes of good reduction \cite[Chapter VII]{Silverman}. Since $\gcd(13,21)=1$, $E(\Q)$ has trivial torsion.\index{torsion}

For the rank, consider the dual $3$-isogenies\index{descent via a 3-isogeny@descent via a $3$-isogeny} $\phi:E\to E'$ and $\psi:E'\to E$, where $E':Y^2=X^3+1998$. Their composites are multiplication by $3$; the geometric kernels consist of the identity and the two points with first coordinate zero on each source curve. Neither kernel has a nonidentity rational point. The classical descent maps for these isogenies \cite[Section 1, Definition 1.3 and Proposition 1.4]{CohenPazuki} are, up to conjugation of the chosen square roots,
\[
 \alpha:E(\Q)\longrightarrow K_1^*/K_1^{*3},\qquad
 \alpha':E'(\Q)\longrightarrow K_2^*/K_2^{*3},
\]
where $K_1=\Q(\sqrt{-74})$ and $K_2=\Q(\sqrt{222})$, and $K_i^*/K_i^{*3}$ is the multiplicative group modulo cubes; they send affine points to the classes of $y+\sqrt{-74}$ and $Y+3\sqrt{222}$, respectively, and the identities to the unit class. Their kernels are $\psi(E'(\Q))$ and $\phi(E(\Q))$, respectively. Since both isogenies are injective on rational points, the index identity gives
\[
 3^{\operatorname{rank}E(\Q)}=[E(\Q):3E(\Q)]=|\operatorname{im}\alpha|\,|\operatorname{im}\alpha'|.
\]
Indeed, $3E(\Q)=\psi\phi(E(\Q))$, and the intermediate subgroup is $\psi(E'(\Q))$. See also \cite[Proposition 2.2]{CohenPazuki}. Only an upper bound on each image is needed.

Write $b=-74$ or $b=1998$, $K=\Q(\sqrt b)$, and let $(x,y)$ be an affine rational point on $y^2=x^3+b$. Neither $b$ is a rational square, so $x\ne0$ and $\theta=y+\sqrt b$ is nonzero. Its norm is $x^3$. For a prime ideal $\mathfrak p$ of the ring of integers of $K$, write $\nu_{\mathfrak p}(\theta)$ for its exponent in the factorization of the principal fractional ideal $(\theta)$. We claim that every such valuation is divisible by $3$. At an inert or ramified rational prime, the norm identity gives $f\nu_{\mathfrak p}(\theta)=3\nu_p(x)$, where $f\in\{1,2\}$ is the residue degree. At a split prime $p$, if $\nu_p(x)<0$, then $2\nu_p(y)=3\nu_p(x)$ and $\nu_p(y)$ is a negative multiple of $3$. The integral element $\sqrt b$ cannot change that valuation at either prime above $p$. If $\nu_p(x)\geq0$, both $y+\sqrt b$ and $y-\sqrt b$ are integral at those primes. No split prime divides $2\sqrt b$: the rational primes dividing its norm are among $2,37$ for $K_1$ and $2,3,37$ for $K_2$, and all are ramified in their respective fields, of discriminants $-296$ and $888$. Hence $\theta$ and its conjugate cannot both have positive valuation at the same split prime. Their product is $x^3$, so again every valuation of $\theta$ is divisible by $3$. This proves $(\theta)=\mathfrak a^3$ for a fractional ideal $\mathfrak a$ of $K$.

Both class groups have trivial $3$-torsion. For $K_1$, the reduced primitive positive definite forms of discriminant $-296$ are
\[
 (1,0,74),\ (2,0,37),\ (3,\pm2,25),\ (5,\pm2,15),\ (6,\pm4,13),\ (9,\pm8,10).
\]
The list follows by checking $1\leq A\leq\sqrt{296/3}$ and $-A\leq B\leq A\leq C$ with $B^2-4AC=-296$, using the usual boundary convention for reduced forms. It gives class number $10$.\index{class number} In the real quadratic field $K_2$, Minkowski's bound\index{Minkowski bound} $\sqrt{888}/2<15$ shows that the ordinary ideal class group is generated by prime ideals of norm at most $14$. The primes $5,7,11$ are inert (their squares exceed $14$); $2,3$ ramify and $13$ splits. If $\mathfrak p_2,\mathfrak p_3$ are the ramified primes, $\mathfrak p_2^2=(2)$ and $\mathfrak p_3^2=(3)$. Moreover $15+\sqrt{222}$ has norm $3$, and $14+\sqrt{222}$ has norm $-26$. The former makes $\mathfrak p_3$ principal; the latter exhibits a prime above $13$ whose class is the inverse of $[\mathfrak p_2]$. The conjugate prime above $13$ has the inverse class. Thus every class is generated by $[\mathfrak p_2]$, of order dividing two. We do not need to decide whether this class is trivial.

Since $[\mathfrak a]^3=1$, we now know $\mathfrak a$ is principal. Writing $\mathfrak a=(\gamma)$ gives $\theta=\epsilon\gamma^3$ for a unit $\epsilon$. The imaginary quadratic field $K_1$ has units $\{1,-1\}$, both cubes; by Dirichlet's unit theorem the units of $K_2$ are generated by $-1$ and one fundamental unit, and there are at most three unit classes modulo cubes. It follows that $|\operatorname{im}\alpha|=1$ and $|\operatorname{im}\alpha'|\leq3$. The rank is consequently at most one. The affine point $\GG=(99,985)$ is rational and, by the torsion calculation, has infinite order, so the rank is exactly one. Its primitivity is established in Section~\ref{sec:height} by the uniform height estimate. This completes the justification of \eqref{eq:group} using the classical descent theorem and elementary arithmetic of the two quadratic fields, without a numerical rank or height computation.

\begingroup
\small
\raggedright

\endgroup
\printindex
\end{document}